\documentclass[11pt]{article}
\usepackage{amssymb}
\usepackage{amsmath}
\usepackage{amsfonts}
\usepackage{mathrsfs}
\usepackage{textcomp}
\usepackage{amsthm}
\usepackage{appendix}
\usepackage{dsfont}
\usepackage{xcolor}
\usepackage{verbatim}
\usepackage{enumitem}
\usepackage{url}
\usepackage{cases}
\usepackage[T1]{fontenc}
\date{}

\newtheorem{Theorem}{Theorem}[section]
\newtheorem{Lemma}{Lemma}[section]

\newtheorem{Corollary}{Corollary}[section]
\newtheorem{Definition}{Definition}[section]
\newtheorem{Proposition}{Proposition}[section]

\numberwithin{equation}{section} \theoremstyle{plain}

\renewcommand{\theequation}{\arabic{section}.\arabic{equation}}

\allowdisplaybreaks[4]

\def\R{{\textbf{R}}}

\def\f{\frac}
\def\R{\mathbb R}
\def\N{\mathbb N}

\def\E{\mathcal E}
\newcommand{\Cc}{\mathcal C}

\title{Uniqueness and Nondegeneracy of Ground States for Critical-Hardy Fractional Equations}

\author{YuJin Guo\thanks{School of Mathematics and Statistics,  Key Laboratory of Nonlinear Analysis $\&$ Applications (Ministry of Education), Central China Normal University, Wuhan 430079, P. R. China. Y. J. Guo is partially supported by
NSF of China (Grants 12225106 and 12371113), and National Key R $\&$ D Program of China (Grant 2023YFA1010001). Email: yguo@ccnu.edu.cn.},
Shuang Wu\thanks{School of Mathematics and Statistics,  and Hubei Key Laboratory of Mathematical Sciences, Central China Normal University, Wuhan 430079, P. R. China. Email: swu@mails.ccnu.edu.cn.},
and  Hongfei Zhang\thanks{School of Mathematics and Statistics, and  Hubei Key Laboratory of Mathematical Sciences, Central China Normal University, Wuhan 430079, P. R. China. H. F. Zhang is partially supported by Doctoral Student Program of the Young S $\&$ T Talents Cultivation Project, CAST. Email: hfzhang@mails.ccnu.edu.cn.} }

\date{\today}

\begin{document}

\maketitle
\begin{abstract}
We study ground states of nonlinear fractional Schr\"odinger equations with the critical Hardy potential in $\R^d$, where $d\geq3$. We prove the existence and qualitative properties of ground states, based on which we analyze the uniqueness and nondegeneracy of ground states. As a byproduct, we also investigate the  optimizers of the critical Hardy--Gagliardo--Nirenberg inequality, from which a sharp Keller--Hardy--Lieb--Thirring inequality is analytically derived.
\end{abstract}


\section{Introduction and main results}

We consider ground states of the following nonlinear Schr\"odinger equation with Hardy potential
\begin{equation}\label{1.6}
(-\Delta)^sQ-\tau\Cc_{d,s}|x|^{-2s}Q+Q=Q^{p-1}\ \ \text{in}\ \ \R^d,
\end{equation}
where $d\geq3$, $0<s\leq1$, $0\leq\tau\leq1$, $2<p<2_s^*:=2d/(d-2s)$, and the constant
\begin{equation*}
\Cc_{d,s}:=2^{2s}\frac{\Gamma((d+2s)/4)^2}{\Gamma((d-2s)/4)^2}>0
\end{equation*}
is the sharp constant (cf. \cite{FLSHardy2008}) of Hardy's inequality
\begin{equation}\label{1.2}
\|(-\Delta)^{s/2}u\|_2^2\geq\Cc_{d,s}\int_{\R^d}|x|^{-2s}|u(x)|^2 dx,\ \ \forall u\in H^s(\R^d).
\end{equation}
Here the parameter $s\in(0,1]$ distinguishes the local and fractional operators, while $\tau\in[0,1]$ measures the strength of the attractive Hardy potential. Specially, the value $\tau=1$ of (\ref{1.6}) is the sharp Hardy coupling, which is referred to as the criticality of (\ref{1.6}) in the present paper.

The study of ground states for \eqref{1.6} is motivated by their connections with spectral optimization problems. More precisely, the Hardy--Lieb--Thirring inequality (cf. \cite{Frank2009,FLSHardy2008}) bounds  the Riesz sum of the negative eigenvalues for the Schr\"odinger  operator $A_{s,1}-V$ in $L^2(\R^d)$, where
\begin{equation}\label{A}
A_{s,\tau}:=(-\Delta)^s-\tau\Cc_{d,s}|x|^{-2s}.
\end{equation}
Particularly, the relation between spectral optimizers of the ordinary operator $A_{1,0}-V$ and ground states of the related equation was discussed in the finite-rank Lieb--Thirring theory of Frank, Gontier and Lewin \cite{FrankGontierLewin2021,FrankGontierLewin2025}, where the spectral optimization problem in the rank-one case is equivalent to a sharp Gagliardo--Nirenberg inequality. In the critical fractional Hardy setting, the corresponding one-particle problem is similarly related to a sharp Hardy--Gagliardo--Nirenberg (HGN) inequality.


When $\tau=0$, equation \eqref{1.6} reduces to the following form
\begin{equation}\label{1.3}
(-\Delta)^s Q+Q=Q^{p-1}\ \ \text{in}\ \ \R^d.
\end{equation}
For $s=1$, Weinstein \cite{Weinstein1983,Weinstein1985} established the connection between ground states of \eqref{1.3}, the sharp Gagliardo--Nirenberg inequality, and the spectrum of the linearized operator for \eqref{1.3}, while the authors in \cite{ChangGustafsonNakanishiTsai2007,Kwong1989} proved the uniqueness and nondegeneracy of   positive radial solutions for \eqref{1.3}. As for $0<s<1$, the authors in \cite{FallValdinoci2014,FrankLenzmann2013,FLS2016} proved the uniqueness and nondegeneracy of ground states for \eqref{1.3}. We remark that the uniqueness theorem of \cite{FLS2016} also gives the complete classification on the optimizers of the corresponding fractional Gagliardo--Nirenberg inequality for \eqref{1.3}.
On the other hand, the Hardy potential $|x|^{-2s}$ arises naturally in quantum mechanics, see \cite{BogdanJakubowskiKimPilarczyk2023,MukherjeeNamNguyen2021}. The potential $|x|^{-2s}$ of (\ref{1.6}) has the same homogeneity as $(-\Delta)^s$, where the sharp Hardy constant determines the threshold for the nonnegativity of the associated quadratic form. Specially, the quadratic form  at the critical coupling $\tau=1$  remains nonnegative, but it does not control any more the homogeneous $H^s$-seminorm of (\ref{1.6}), and the Friedrichs form domain  of (\ref{1.6}) strictly contains $H^s(\mathbb R^d)$, see \cite[Section~2.1]{FLSHardy2008}. This loss of coercivity presents therefore a core difficulty in the analysis of the critical Hardy problem (\ref{1.6}) with $\tau =1$.

When $s=1$ and $0<\tau\leq1$, equation \eqref{1.6} becomes the following local Hardy--Schr\"odinger equation
\begin{equation}\label{1.4}
\left(-\Delta-\tau\frac{(d-2)^2}{4}|x|^{-2}\right)Q+Q
=Q^{p-1}\ \,\ \text{in}\ \ \R^d.
\end{equation}
For $\tau=1$ and $2<p<2+4/d$, the existence and singularity of positive radial solutions for \eqref{1.4} were analyzed in  \cite{TrachanasZographopoulos2015}. While for $\tau=1$ and $2<p<2d/(d-2)$, the authors in \cite[Theorems~1 and~2]{MukherjeeNamNguyen2021} proved the uniqueness of   positive radial solutions for \eqref{1.4} and classified all optimizers of the associated sharp HGN inequality, which were further extended in \cite{MukherjeeNamNguyen2021,YangZengZhang2026} to the non-critical case $0<\tau<1$.

When $0<s<1$ and $0<\tau\leq1$,
equation \eqref{1.6} exhibits  simultaneously both the nonlocality and a Hardy-type singularity, which naturally raises the question on whether the existence, uniqueness and other analytical properties of ground states for \eqref{1.6} still hold in this fractional Hardy setting, particularly at the critical case $\tau=1$.
The main purpose of the present paper is to address this issue. As a byproduct, we shall also investigate the optimal constant and classify all optimizers of the critical Hardy--Gagliardo--Nirenberg inequality, from which a sharp Keller--Hardy--Lieb--Thirring inequality is finally derived.

\subsection{Main results}

The main purpose of this subsection is to introduce the main results of the present paper. From now on, we define for $0<s<1$ and $0\leq\tau\leq1$,
\begin{align}
\mathfrak q_{s,\tau}[u,v]&:=\int_{\R^d}\Big[(-\Delta)^{s/2}u (-\Delta)^{s/2}v-\tau\Cc_{d,s}|x|^{-2s}uv\Big] dx,\label{2.1}\\
\mathfrak h_{s,\tau}[u,v]&:=\mathfrak q_{s,\tau}[u,v]+\int_{\R^d}uvdx , \ \ \, \forall\,u, v\in C_c^\infty(\R^d\setminus\{0\}) .\label{2.2}
\end{align}
For simplicity, denote
$$\mathfrak q_{s,\tau}[u]:=\mathfrak q_{s,\tau}[u, u],\quad \mathfrak h_{s,\tau}[u]:=\mathfrak h_{s,\tau}[u, u].$$

When $0\leq\tau<1$, the fractional Hardy inequality \eqref{1.2} implies that $\mathfrak q_{s,\tau}$ and $\mathfrak h_{s,\tau}$ are also well defined on $H^s(\R^d)$, and their closures have the form domain $H^s(\R^d)$. Moreover, $\mathfrak h_{s,\tau}^{1/2}$ defines a norm equivalent to the standard $H^s$-norm, while $\mathfrak q_{s,\tau}$ is comparable to the square of the homogeneous $H^s$-seminorm. While for $\tau=1$, we understand $\mathfrak q_{s,1}$ as the Friedrichs closure of the form in \eqref{2.1}, and its form domain, equivalently the form domain of
$\mathfrak h_{s,1}$, can be denoted by
\begin{equation}\label{eq:Qs}
\mathcal Q^s:=\Big\{u\in L^2(\R^d):~~\int_{\R^d}\big(|(-\Delta)^{s/2}u|^2-\Cc_{d,s}|x|^{-2s}|u|^2\big)dx<\infty\Big\},
\end{equation}
which is equipped with the inner product
\begin{equation*}
\langle u,v\rangle_{\mathcal Q^s}:=\int_{\R^d}\big[(-\Delta)^{s/2}u(-\Delta)^{s/2}v-\Cc_{d,s}|x|^{-2s}uv+uv\big]dx
\end{equation*}
and the norm $\|u\|_{\mathcal Q^s}:=\langle u,u\rangle_{\mathcal Q^s}^{1/2}=\mathfrak h_{s,1}[u]^{1/2}$. Then $\mathcal Q^s$ is a completion (cf.\cite{FLSHardy2008}) of $C_c^\infty(\R^d\setminus\{0\})$ under the norm $\|\cdot\|_{\mathcal Q^s}$.
Consequently, the form domain of $\mathfrak h_{s,\tau}$ satisfies
\begin{equation}\label{0D}
D(\mathfrak h_{s,\tau})=\begin{cases}
H^s(\R^d),  &0\leq\tau<1;\\
\mathcal Q^s,&\tau=1.
\end{cases}
\end{equation}
We shall study ground state  of \eqref{1.6} in the following sense.

\begin{Definition} \label{def:morse-one}
For $\tau\in[0,1]$, assume that $Q_\tau\in D(\mathfrak h_{s,\tau})$ is a nonnegative and nontrivial solution of equation~\eqref{1.6}, and let
\begin{equation}\label{1.8}
L_{+,\tau}:=(-\Delta)^s-\tau\Cc_{d,s}|x|^{-2s}+1-(p-1)Q_{\tau}^{p-2}
\end{equation}
be the associated self-adjoint linearized operator in $L^2(\R^d)$.  The solution $Q_\tau$ is then called a ground state of equation \eqref{1.6}, if the operator $L_{+,\tau}$ in $L^2(\R^d)$ has Morse index $1$, i.e., $L_{+,\tau}$ in $L^2(\R^d)$ has exactly one
negative eigenvalue (counting multiplicity).
\end{Definition}

One can note from the strong maximum principle (cf. \cite[Theorem 1.1]{JarohsWeth2019}) that any nonnegative and nontrivial solution $Q_\tau\in D(\mathfrak h_{s,\tau})$ of \eqref{1.6} must be positive in $\R^d$. Following Definition~\ref{def:morse-one}, ground states of \eqref{1.6} are essentially positive for any $0<\tau\leq1$.
We first give the following existence and qualitative properties of ground states for \eqref{1.6} with $0<\tau\leq1$.

\begin{Theorem}\label{thm1.1}
Suppose  $0<s<1$, $0<\tau\leq1$ and $2<p<2_s^*:=\frac{2d}{d-2s}$, where $d\geq3$. Then \eqref{1.6} admits a ground state $Q_\tau\in D(\mathfrak h_{s,\tau})$.
Moreover, any ground state $Q_\tau\in D(\mathfrak h_{s,\tau})$ of \eqref{1.6} must be radially symmetric about the origin, strictly decreasing in $r=|x|$, and satisfies
\begin{equation}\label{1.11}
\begin{aligned}
Q_\tau\in C^\infty(\R^d\setminus\{0\}),\ \ \lim_{|x|\to0}|x|^{(d-2s)/2}Q_1(x):=\ell_Q\in(0,\infty).
\end{aligned}
\end{equation}
\end{Theorem}

  To prove the existence of Theorem \ref{thm1.1}, we define
\begin{equation}\label{1.7}
\nu_{s,\tau}:=\inf\left\{\mathfrak h_{s,\tau}[u]:\, u\in D(\mathfrak h_{s,\tau}),\ \, \|u\|_p=1\right\},
\end{equation}
where $\mathfrak h_{s,\tau}[u]:=\mathfrak h_{s,\tau}[u, u]$ is as in (\ref{2.2}).
The existence and regularity of ground states for \eqref{1.6} then follow from the standard variational and bootstrap arguments (cf. \cite{Frank2009,FLS2016}). As for the symmetry of Theorem \ref{thm1.1}, unfortunately the argument of \cite{BhaktaGangulyMontoro2023} is not applicable to the case $\tau=1$. For this reason, we shall address the symmetry  of Theorem \ref{thm1.1} by applying the pole lower estimate of \cite{BogdanEtAl2019}, together with a pointwise antisymmetric argument. Moreover,  the precise pole asymptotics of \eqref{1.11} is proved from the normalized heat-kernel limit of \cite[Theorem~1.1]{BogdanJakubowskiKimPilarczyk2023} and the nonlinear resolvent representation
of ground states for \eqref{1.6}. We also note from \eqref{1.2} and \eqref{1.11} that $Q_1\notin H^s(\R^d)$. The proof of Theorem~\ref{thm1.1} is left to Appendix~\ref{sec:forms} for the reader's convenience.

The main purpose of the present paper is concerned with the following uniqueness and nondegeneracy of ground state  for \eqref{1.6}.

\begin{Theorem}
\label{thm:main}
Suppose  $0<s<1$, $0<\tau\leq1$ and $2<p<2_s^*:=\frac{2d}{d-2s}$, where $d\geq3$. Then  \eqref{1.6}  has a unique  ground state $Q_\tau$, and the Friedrichs extension of the linearized operator $L_{+,\tau}$ around $Q_\tau$ in $L^2(\R^d)$ satisfies
\begin{equation}\label{1.14}
\ker L_{+,\tau}=\{0\}.
\end{equation}
\end{Theorem}

We remark that when $s=1$, Mukherjee, Nam and Nguyen proved in \cite{MukherjeeNamNguyen2021} the uniqueness of positive radial solutions for \eqref{1.6} in $\mathcal Q^1$. We now sketch the proof of Theorem \ref{thm:main} by the following three steps.

In the first step of proving Theorem \ref{thm:main}, we establish the radial nondegeneracy of ground states $Q_\tau$ for \eqref{1.6}, where $0<\tau\leq1$. Unfortunately, the radial nodal theorem in \cite{FLS2016} is not applicable, due to the appearance of the singular Hardy potential. To overcome this difficulty, we approximate the singular Hardy potential by bounded radial monotone potentials. Applying a compactness argument, we then employ \cite[Theorem~2.3]{FLS2016} to show that any nontrivial radial kernel function of $L_{+,\tau}$ in $L^2(\R^d)$ must change exactly one sign. A further difficulty however arises from the critical coupling $\tau=1$, where the kinetic and Hardy energies of $Q_1$ are separately infinite. We shall resolve this issue by introducing a bilinear argument in logarithmic radial coordinates to establish the following orthogonality identity
\begin{equation*}
\int_{\R^d}Q_\tau Z\,dx=0,\ \ \ Z\in \ker L_{+,\tau}\big|_{\mathrm{rad}},\ \ \ \tau\in(0,1].
\end{equation*}
Since $\int_{\R^d}Q_\tau^{p-1}Z\,dx=0$,
the strict radial monotonicity of $Q_\tau$ then yields a contradiction, which therefore proves in Proposition~\ref{prop:radialnon} the radial nondegeneracy of ground states $Q_\tau$ for \eqref{1.6}.

As the second step  of proving Theorem \ref{thm:main}, we exclude all nonradial zero modes of $L_{+,\tau}$ in $L^2(\R^d)$ for any $0<\tau\leq1$. Differentiating \eqref{1.6} away from the origin, it actually yields that
\begin{equation*}
L_{+,\tau}(-\partial_jQ_\tau)
=2s\tau\Cc_{d,s}|x|^{-2s-2}Q_\tau(x)x_j \ \,\ \text{in}\ \, \R^d\setminus\{0\},
\quad j=1,\cdots,d.
\end{equation*}
We can derive from above that $-Q_\tau'(r)/r>0$ holds for all $r>0$, which can further yield the strict positivity of the quadratic form for $L_{+,\tau}$ in the first angular momentum channel. However, when $\tau=1$, the corresponding form domain of $L_{+,\tau}$ is $\mathcal Q^s$, instead of $H^s(\R^d)$, which yields the difficulty that the angular comparison of \cite{FLS2016} is not applicable. To handle this difficulty, we shall follow the improved angular Hardy inequalities of \cite{BogdanMerz2024} to identify the nonradial form domains with the corresponding subspaces of $H^s(\R^d)$. The angular comparison of \cite[Lemma~C.5]{FLS2016} then excludes   zero modes of $L_{+,\tau}$ in all higher angular momentum channels. Together with the radial nondegeneracy in the first step, this finally proves that $\ker L_{+,\tau}=\{0\}$ holds for $0<\tau\leq1$, $i.e.$, which hence gives the nondegeneracy of Theorem \ref{thm:main}.

The third step  of proving Theorem \ref{thm:main} is to prove the uniqueness of ground states $Q_\tau$ for \eqref{1.6} by a continuation argument, where $0<\tau\leq 1$. Unfortunately, since the kinetic and Hardy energies of $Q_1$ are separately infinite, and $\mathcal Q^s$ depends on $s$, the continuation method of Frank, Lenzmann and Silvestre \cite{FLS2016}, which varies $s\in(0,1]$, cannot be applied directly  at $\tau=1$. To overcome this issue, we  fix $s\in(0,1)$ and vary $\tau\in[0,1]$.
The new difficulty however appears that $D(\mathfrak h_{s,\tau})$ changes at $\tau=1$. For this reason, we shall employ Mosco's convergence in \cite{Mosco1994} and the radial compactness of the
embedding $\mathcal Q^s_{\rm rad}(\R^d)\hookrightarrow L^p(\R^d)~(2<p<2_s^*)$ to derive the operator-norm continuity of the resolvent
$$R_\tau:=\big((-\Delta)^s-\tau\Cc_{d,s}|x|^{-2s}+1\big)^{-1},\ \ \mbox{where} \ \,\tau\in[0,1].$$
Applying the nondegeneracy of $Q_\tau$, this then yields locally a unique branch of each ground state $Q_\tau$ for \eqref{1.6} near $\tau=1$. We further use a priori estimates and the strong compactness
of ground-state branches away from the critical coupling $\tau=1$ to extend these unique branches  throughout $\tau\in[0, 1]$, along each of which the Morse indexes of ground states remain one. In view of \cite[Theorem~3.4]{FLS2016}, all branches obtained above must coincide at $\tau=0$. The local uniqueness of all branches at $\tau=0$ finally yields the uniqueness of ground states for \eqref{1.6} with any $0<\tau \leq 1$. We leave the detailed proof of Theorem~\ref{thm:main} to Sections \ref{sec:nondegeneracy} and \ref{sec:cont-fixed-order}.

We next focus on the applications of Theorem \ref{thm:main}. Define
\begin{equation}\label{01}
\theta:=\frac{d(p-2)}{2sp}\in(0,1),
\end{equation}
where $0<s<1$,   $2<p<2_s^*:=\frac{2d}{d-2s}$ and $d\geq3$ are as in \eqref{1.6}.
Applying Theorem \ref{thm:main}, we are able to analyze the optimizers of the following Hardy--Gagliardo--Nirenberg (HGN) inequality.

\begin{Theorem}\label{thm:sharp-HGN}
Suppose  $0<s<1$, $0<\tau\leq1$ and $2<p<2_s^*:=\frac{2d}{d-2s}$, where $d\geq3$. Consider the Hardy--Gagliardo--Nirenberg inequality
\begin{equation}\label{1.17}
\|u\|_p^2\leq C_{\rm HGN}\Big(\int_{\R^d}\Big[\big|(-\Delta)^{s/2}u\big|^2-\frac{\Cc_{d,s}}{|x|^{2s}}|u|^2\Big] dx\Big)^\theta\|u\|_2^{2(1-\theta)},\ \ \forall~u\in\mathcal{Q}^s,
\end{equation}
where $\mathcal{Q}^s$ is defined by \eqref{eq:Qs}, and let $Q_1>0$ be the unique ground state of \eqref{1.6} with $\tau=1$. Then we have
\begin{enumerate}
\item[(1).]
The best constant $C_{\rm HGN}>0$ satisfies
\begin{equation}\label{eq:sharp-HGN-constant-nu}
C_{\rm HGN}=\frac{1}{\|Q_1\|_p^{p-2}\theta^\theta(1-\theta)^{1-\theta}}.
\end{equation}

\item[(2).] The inequality \eqref{1.17} is attained   by
\begin{equation}\label{1.21}
u(x)\equiv \pm a\lambda^{d/2}Q_1(\lambda x)\ \ \text{in \ } \R^d\setminus\{0\},\ \ \forall\,a>0,~\lambda>0.
\end{equation}
\end{enumerate}
\end{Theorem}

Theorem \ref{thm:sharp-HGN} establishes a sharp one-particle interpolation
inequality, and we refer to \cite{DuongLeNamNguyen2025}
for the related many-particle version of \eqref{1.17}.
It is worth emphasizing that the inequality of the form \eqref{1.17} was established in
\cite[Corollary~2.5 and (2.8)]{FLSHardy2008}, whose sharp
constant and optimizers were however not determined.
Our result identifies the sharp constant and characterizes
all nonzero optimizers of \eqref{1.17}.
Since the sharp inequality \eqref{1.17} follows from the variational problem \eqref{1.7} by a scaling argument, we prove Theorem \ref{thm:sharp-HGN} (1) in Appendix~\ref{sec:forms}. To classify all nonzero optimizers of \eqref{1.17}, we follow the argument of \cite{MukherjeeNamNguyen2021,Weinstein1983}, for which the difficulty arises from the critical Hardy potential. We shall address this issue by using the ground-state
representation of \cite{FrankSeiringer2008}. For any nonzero optimizer $u$, the optimality implies that $\mathfrak q_{s,1}[|u|]=\mathfrak q_{s,1}[u]$, where $\mathfrak q_{s,1}$ is defined by \eqref{2.1}. After a suitable rescaling, $|u|$ becomes a positive ground state of \eqref{1.6} with $\tau=1$. By the uniqueness result of Theorem~\ref{thm:main}, this ground state
must coincide with $Q_1$. This thus proves Theorem \ref{thm:sharp-HGN} (2).

As an application of Theorem \ref{thm:sharp-HGN}, we finally derive a sharp spectrum inequality. Towards this aim, we denote $\lambda_1(H):=\inf\sigma(H)$ the bottom of the spectrum for a self-adjoint operator $H$ bounded from below, and define $t_\pm:=\max\{\pm t,0\}$ for $t\in\R$.

\begin{Corollary} \label{cor:one-particle-HLT}
Suppose  $0<s<1$, $0<\tau\leq1$ and $d\geq3$, and let  $Q_1>0$ be the unique ground state of \eqref{1.6} with $\tau=1$. Then for every $\gamma>0$,
\begin{equation}\label{eq:KHLT}
\Big[\lambda_1\big((-\Delta)^s-\Cc_{d,s}|x|^{-2s}-V\big)\Big]_-^\gamma\leq \nu_{s,1}^{-(\gamma+\frac{d}{2s})}\int_{\R^d}V_+(x)^{\gamma+\frac{d}{2s}} dx,~~~\forall~V(x)\in L^{\gamma+\frac{d}{2s}}(\R^d),
\end{equation}
where the constant $\nu_{s,1}^{-(\gamma+\frac{d}{2s})}>0$ is optimal and  $\nu_{s,1}>0$ is defined by \eqref{1.7}. Moreover, \eqref{eq:KHLT} is attained precisely  by
\begin{equation}\label{eq:one-particle-HLT-optimizers}
V(x)\equiv \lambda^{2s}Q_1(\lambda x)^{p-2}\text{ \ in \  } \R^d\setminus\{0\},\ \ \forall\,\lambda>0,
\end{equation}
where $p=\frac{2(d+2s\gamma)}{d+2s\gamma-2s}\in(2,2_s^*)$.
\end{Corollary}

Corollary~\ref{cor:one-particle-HLT} establishes a sharp spectrum inequality \eqref{eq:KHLT}. The one-particle optimization problem for the Schr\"odinger operator $-\Delta-V$ goes back to Keller \cite{Keller1961}. The proof of Corollary~\ref{cor:one-particle-HLT} follows from a similar argument of \cite{CarlenFrankLieb2014,FrankGontierLewin2021}. More precisely, we combine the sharp HGN inequality \eqref{1.17} with a scalar optimization to derive the sharp
spectral bound of \eqref{eq:KHLT}, and then use Theorem~\ref{thm:sharp-HGN} to characterize all optimizing potentials.

The paper is organized as follows.
In Section~\ref{sec:nondegeneracy}, we establish the nondegeneracy of ground states for \eqref{1.6}.
In Section~\ref{sec:cont-fixed-order}, we prove the uniqueness of ground states for  \eqref{1.6} by the continuation in the Hardy coupling at fixed $s$, thereby completing the proof of Theorem~\ref{thm:main} in Subsection~\ref{sub3.1}.  Subsection~\ref{sec:completion} gives the optimizer classification of Theorem~\ref{thm:sharp-HGN} and proves Corollary~\ref{cor:one-particle-HLT}. Appendix~\ref{sec:forms} proves Theorem~\ref{thm1.1} and establishes the sharp interpolation inequality \eqref{1.17}.

\section{Nondegeneracy of Ground States}\label{sec:nondegeneracy}
In this section, we mainly prove Theorem \ref{thm:fullkernel} on the nondegeneracy of the linearized operator $L_{+,\tau}$ at every ground state $Q:=Q_\tau$ of \eqref{1.6}, where $0<\tau\leq1$.

Let $Z\in D(\mathfrak h_{s,\tau})$ be a radial weak solution of the equation
\begin{equation} \label{eq:radker}
 L_{+,\tau}Z:=(-\Delta)^sZ-\tau\Cc_{d,s}|x|^{-2s}Z+Z-(p-1)Q^{p-2}Z=0\ \ \text{in}\ \ \R^d.
\end{equation}
One can check that $Z\in C^\infty(\R^d\setminus\{0\})$ by a similar argument of Lemma \ref{lem3.1}. It then follows from  \eqref{1.6} and \eqref{eq:radker} that
\begin{equation}\label{eq:1moment}
\int_{\R^d}Q^{p-1}Z dx=0.
\end{equation}
Define for $f,g\in C_c^\infty(\R)$,
\begin{equation*}
U(x):=|x|^{-(d-2s)/2}f(-\log |x|),\ \ V(x):=|x|^{-(d-2s)/2}g(-\log |x|).
\end{equation*}
Applying \cite[Theorem~1.2]{FrankSeiringer2008} with $p=2$ and using polar coordinates, we obtain that
\begin{equation}\label{eq:cyform}
\begin{aligned}
\frac{\mathfrak q_{s,1}[U,V]}{|\mathbb S^{d-1}|}
&=\mathfrak j_s[f,g]\\
&:=\frac{c_{d,s}}2\iint_{\R^2}\big[f(t)-f(u)\big]\big[g(t)-g(u)\big]K_{d,s}(t-u) dt du,
\end{aligned}
\end{equation}
where $d\omega$ denotes the standard surface measure on $\mathbb S^{d-1}$, and
\begin{equation*}
K_{d,s}(\xi):=\int_{\mathbb S^{d-1}}\frac{d\omega}{\{2(\cosh\xi-\omega_1)\}^{(d+2s)/2}}.
\end{equation*}
Denote $\mathfrak j_s[f]:=\mathfrak j_s[f,f]$. By the density theory, \eqref{eq:cyform} can be extended to radial functions $U,V\in D(\mathfrak h_{s,\tau})$.

Recall from \cite[Proposition~2.1]{AoEtAl2020} that
\begin{equation}\label{eq:kernound1}
K_{d,s}(\xi)\sim |\xi|^{-1-2s},
~~~0<|\xi|\leq1,
\end{equation}
where and below $\sim$ means that the two sides are comparable
up to positive constants depending only on $d$ and $s$.
Differentiating the formula of $K_{d,s}(\xi)$ yields that
\begin{equation}\label{eq:kernelbound2}
\frac{|K_{d,s}'(\xi)|}{K_{d,s}(\xi)}\leq C\bigl(1+|\xi|^{-1}\bigr),~~~\xi\neq0.
\end{equation}
Set for $0\leq\tau\leq1$,
\begin{equation}\label{eq:cytauform}
\mathfrak t_{s,\tau}[f,g]:=\mathfrak j_s[f,g]+(1-\tau)\Cc_{d,s}\int_\R fg dt.
\end{equation}
Set $r=e^{-t}$, and define
\begin{equation*}
q(t):=r^{(d-2s)/2}Q(r),\ \ \zeta(t):=r^{(d-2s)/2}Z(r),\ \  \beta:=d-p(d-2s)/2>0.
\end{equation*}
Multiplying \eqref{1.6} and \eqref{eq:radker} by $|x|^{-(d-2s)/2}\phi(-\log |x|)$,  one can calculate from \eqref{eq:cytauform} that for every $\phi\in C_c^\infty(\R)$,
\begin{equation}\label{eq:cyQ}
\mathfrak t_{s,\tau}[q,\phi]+\int_\R e^{-2st}q\phi dt=\int_\R e^{-\beta t}q^{p-1}\phi dt,
\end{equation}
 and
 \begin{equation}\label{eq:cyz}
\mathfrak t_{s,\tau}[\zeta,\phi]+\int_\R e^{-2st}\zeta\phi dt=(p-1)\int_\R e^{-\beta t}q^{p-2}\zeta\phi dt.
\end{equation}
Let $\chi\in C_c^\infty(\R,[0,1])$ satisfy $\chi\equiv1$ on $[-1,1]$, and $\chi\equiv0$ outside $[-2,2]$. Define $\chi_R(t):=\chi(t/R)$ for $R\geq1$.  Taking $\phi=\chi_R\zeta'$ in \eqref{eq:cyQ} and $\phi=\chi_Rq'$ in \eqref{eq:cyz}, we then deduce that
\begin{equation}\label{eq:fullremainder}
\begin{aligned}
&2s\int_\R\chi_Re^{-2st}q\zeta dt-\beta\int_\R\chi_Re^{-\beta t}q^{p-1}\zeta dt\\
=&-\mathfrak j_s[q,\chi_R\zeta']-\mathfrak j_s[\zeta,\chi_Rq']+(1-\tau)\Cc_{d,s}\int_\R\chi_R'q\zeta dt\\
&+\int_\R\chi_R'e^{-2st}q\zeta dt-\int_\R\chi_R'e^{-\beta t}q^{p-1}\zeta dt.
\end{aligned}
\end{equation}

Since $\mathfrak j_s[q]+\mathfrak j_s[\zeta]<\infty$, the Cauchy--Schwarz inequality yields that
\begin{equation*}
\begin{aligned}
\iint_{\R^2}K_{d,s}(t-u)|q(t)-q(u)| |\zeta(t)-\zeta(u)| dt du
\leq\frac{2}{c_{d,s}}\mathfrak j_s[q]^{1/2}\mathfrak j_s[\zeta]^{1/2}<\infty.
\end{aligned}
\end{equation*}
Choose an even function $\varphi\in C^\infty(\R,[0,1])$ satisfying $\varphi\equiv0$ on $[-1,1]$ and $\varphi\equiv1$ outside $[-2,2]$. For any fixed $R>1$, it yields that
\begin{equation*}
\begin{aligned}
&\frac{c_{d,s}}2\iint_{\R^2}\varphi\left(\frac{t-u}{\varepsilon}\right)K_{d,s}(t-u)\times\Bigl[\big(q(t)-q(u)\big)\{\chi_R(t)\zeta'(t) -\chi_R(u)\zeta'(u)\}\\
&+(\zeta(t)-\zeta(u))\{\chi_R(t)q'(t)-\chi_R(u)q'(u)\}\Bigr] dt du\\
=&-\frac{c_{d,s}}2\iint_{\R^2}\bigg\{\partial_t\Big[\chi_R(t)\varphi\left(\frac{t-u}{\varepsilon}\right)K_{d,s}(t-u)\Big]+\partial_u\Big[\chi_R(u)\varphi\left(\frac{t-u}{\varepsilon}\right)K_{d,s}(t-u)\Big]\bigg\}\\
&\times[q(t)-q(u)]\times[\zeta(t)-\zeta(u)] dt du.
\end{aligned}
\end{equation*}
Since $|\chi_R(t)-\chi_R(u)|\leq C|t-u|/R$, one can check that the absolute value of the above terms containing $\varphi'$ is bounded by
\begin{equation*}
\begin{aligned}
\frac{C}{R}\iint_{\varepsilon\leq|t-u|\leq2\varepsilon}K_{d,s}(t-u)|q(t)-q(u)| |\zeta(t)-\zeta(u)| dt du\to0 \text{ \ as \  }\varepsilon\searrow0.
\end{aligned}
\end{equation*}
Setting $\varepsilon\searrow0$ and using dominated convergence, it yields from \eqref{eq:kernelbound2} that
\begin{equation*}
\begin{aligned}
&\mathfrak j_s[q,\chi_R\zeta']+\mathfrak j_s[\zeta,\chi_Rq']\\
=&-\frac{c_{d,s}}2\iint_{\R^2}K_{d,s}(t-u)\{\chi_R'(t)+\chi_R'(u)\}\times[q(t)-q(u)]\times[\zeta(t)-\zeta(u)] dt du\\
&-\frac{c_{d,s}}2\iint_{\R^2}K_{d,s}'(t-u)\{\chi_R(t)-\chi_R(u)\}\times[q(t)-q(u)]\times[\zeta(t)-\zeta(u)] dt du.
\end{aligned}
\end{equation*}
Using $\|\chi_R'\|_\infty\leq C/R$, we then have
\begin{equation*}
\begin{aligned}
&\big|\mathfrak j_s[q,\chi_R\zeta']+\mathfrak j_s[\zeta,\chi_Rq']\big| \\
\leq& \frac{C}{R}\mathfrak j_s[q]^{1/2}\mathfrak j_s[\zeta]^{1/2}
 +\frac{c_{d,s}}2\iint_{\R^2}|K_{d,s}'(t-u)| |\chi_R(t)-\chi_R(u)| |q(t)-q(u)| |\zeta(t)-\zeta(u)| dt du.
\end{aligned}
\end{equation*}
Since $|\chi_R(t)-\chi_R(u)|\leq C\min\{1,|t-u|/R\}$, it yields from \eqref{eq:kernelbound2} that for $R>1$ and $t\neq u$,
\begin{equation*}
\begin{aligned}
|K_{d,s}'(t-u)| |\chi_R(t)-\chi_R(u)|
\leq&C K_{d,s}(t-u)(1+|t-u|^{-1})\min\{1,|t-u|/R\}\\
\leq&C K_{d,s}(t-u).
\end{aligned}
\end{equation*}
Since $C K_{d,s}(t-u)|q(t)-q(u)| |\zeta(t)-\zeta(u)|$ is integrable on $\R\times\R$, the dominated convergence theorem gives that
\begin{equation*}
\iint_{\R^2}|K_{d,s}'(t-u)| \times|\chi_R(t)-\chi_R(u)|\times|q(t)-q(u)|\times |\zeta(t)-\zeta(u)| dt du\to0\text{ \ as \ }R\to\infty.
\end{equation*}
Consequently, we have
\begin{equation}\label{eq:remainder0}
\mathfrak j_s[q,\chi_R\zeta']+\mathfrak j_s[\zeta,\chi_Rq']\to0 \text{ \ as \  }R\to\infty.
\end{equation}

Setting $R\to\infty$ into \eqref{eq:fullremainder}, we thus derive from \eqref{eq:remainder0} that
\begin{equation*}
2s\int_\R e^{-2st}q\zeta dt=\beta\int_\R e^{-\beta t}q^{p-1}\zeta dt.
\end{equation*}
Using $\beta=d-p(d-2s)/2$ and
\begin{equation*}
|\mathbb S^{d-1}|\int_\R e^{-2st}q\zeta dt=\int_{\R^d}QZ dx,\ \
|\mathbb S^{d-1}|\int_\R e^{-\beta t}q^{p-1}\zeta dt=\int_{\R^d}Q^{p-1}Zdx,
\end{equation*}
 we obtain that
\begin{equation*}
2s\int_{\R^d}QZ dx=\left(d-\frac{p(d-2s)}2\right)\int_{\R^d}Q^{p-1}Z dx.
\end{equation*}
Together with \eqref{eq:1moment}, this yields that
\begin{equation}\label{eq:2moment}
\int_{\R^d}QZ dx=0,
\end{equation}
where $Z\in D(\mathfrak h_{s,\tau})$ is a radial weak solution of \eqref{eq:radker}.

Using orthogonality conditions \eqref{eq:1moment} and  \eqref{eq:2moment}, we now prove the following radial nondegeneracy.

\begin{Proposition} \label{prop:radialnon}
Suppose $0<s<1$, $0\leq\tau\leq1$, and $2<p<2_s^*$. If $Q>0$ is a radial weak solution of \eqref{1.6} and the Morse index of $L_{+,\tau}$ is one, then
\begin{equation*}
\ker L_{+,\tau}\big|_{\mathrm{rad}}=\{0\}.
\end{equation*}
\end{Proposition}

\noindent\textbf{Proof.}
Write the radial linearized operator of  $L_{+,\tau}$ as
\begin{equation*}
H:= L_{+,\tau}|_{rad}=(-\Delta)^s+1-W,\ \ W(x):=\tau\Cc_{d,s}|x|^{-2s}+(p-1)Q(x)^{p-2}.
\end{equation*}
By H\"older's inequality and the HGN inequality \eqref{1.17}, we have for every $u\in D(\mathfrak h_{s,\tau})$,
\begin{equation*}
\int_{\R^d}Q^{p-2}|u|^2 dx\leq\|Q\|_{L^p}^{p-2}\|u\|_{L^p}^2\leq C\mathfrak q_{s,\tau}[u]^\theta
\|u\|_2^{2(1-\theta)},\ \ \theta:=\frac{d(p-2)}{2sp}\in(0,1).
\end{equation*}
Thus, it yields from Young's inequality that for every $\varepsilon>0$,
\begin{equation}\label{eq:nonforms}
\int_{\R^d}Q^{p-2}|u|^2 dx\leq\varepsilon\mathfrak q_{s,\tau}[u]+C_\varepsilon\|u\|_2^2.
\end{equation}
Consequently, by the KLMN theorem (cf. \cite[Theorem 3.19]{Lewin}), the bilinear form
\begin{equation*}
h[u,v]:=\mathfrak h_{s,\tau}[u,v]-(p-1)\int_{\R^d}Q^{p-2}uv dx
\end{equation*}
is closed and lower semibounded on
\begin{equation*}
D(h)=
\begin{cases}
H^s_{\rm rad}(\R^d),&0\leq\tau<1,\\
\mathcal Q^s_{\rm rad}(\R^d),&\tau=1,
\end{cases}
\end{equation*}
whose associated self-adjoint operator is $H$. Here and below, denote $h[u]:=h[u,u]$. Moreover, since $h$ is lower semibounded, we choose $C_0>0$ such that
\begin{equation*}
h[u]\geq -C_0\|u\|_2^2,\quad u\in D(h),
\end{equation*}
where $D(h)$ is equipped with the form norm
\begin{equation*}
\|u\|_{D(h)}^2:=h[u]+(C_0+1)\|u\|_2^2.
\end{equation*}

Since $Q$ is radial and $h[Q]=-(p-2)\|Q\|_p^p<0$, the Morse-index assumption implies that
$\operatorname{ind}_{\rm rad}(H)=1$, where $\operatorname{ind}_{\rm rad}(H)$ denotes the Morse index of $H$ in $L^2_{\rm rad}(\R^d)$. Assume that $0\not\equiv z\in\ker H$. Denote
\begin{equation*}
\lambda_2^{\rm rad}(H):=\inf_{\substack{E\subset D(h)\\ \dim E=2}}\sup_{0\neq u\in E}\frac{h[u]}{\|u\|_2^2}.
\end{equation*}
The Morse-index assumption then implies that $\lambda_2^{\rm rad}(H)\geq0$. Moreover, there exists $v\in D(h)$ such that $h[v]<0$. Since $z\in\ker H$, we have $h[v,z]=0$ and  $h[z]=0$.
Thus, we get that
\begin{equation*}
\sup_{0\neq u\in E}\frac{h[u]}{\|u\|_2^2}=0,\ \ E:=\operatorname{span}\{v,z\},
\end{equation*}
This implies that $\lambda_2^{\rm rad}(H)=0$.
Choose a nondecreasing function
$\Theta\in C^\infty([0,\infty))$ such that
\begin{equation*}
0\leq\Theta'(t)\leq1,\ \ \Theta(t)=t\quad(0\leq t\leq1),\ \ \Theta(t)=2\quad(t\geq3),
\end{equation*}
and define
\begin{equation}\label{eq:smoothcut}
W_n(x):=n\eta\big(|x|/n\big)\Theta\bigl(W(x)/n\bigr),
\end{equation}
where $\eta\in C_c^\infty([0,\infty),[0,1])$ is nonincreasing and satisfies  $\eta\equiv1$ on $[0,1]$ and $\eta\equiv0$ on $[2,\infty)$.
Then
\begin{equation*}
0\leq W_n\leq W,
\ \
W_n=W\text{ \ on \ }\{W\leq n,\ |x|\leq n\},
\ \
W_n(x)\to0\quad \text{as}\quad |x|\to\infty,
\end{equation*}
and $W_n$ is bounded and radial decreasing in $|x|$.

Define
\begin{equation*}
h_n[u]
:=
\|(-\Delta)^{s/2}u\|_2^2+\|u\|_2^2
-\int_{\R^d}W_n|u|^2 dx,
\ \
D(h_n)=H^s_{\rm rad}(\R^d),
\end{equation*}
and let $H_n$ be the associated self-adjoint operator.
Since $W_n\leq W$, we have $h_n[u]\geq h[u]$ for every $u\in H^s_{\rm rad}(\R^d)$. Hence, by the min--max principle and $H^s_{\rm rad}(\R^d)\subset D(h)$, we deduce that
\begin{equation*}
\begin{aligned}
\lambda_{2,n}&:=\lambda_2^{\rm rad}(H_n)=\inf_{\substack{E\subset H^s_{\rm rad}(\R^d)\\ \dim E=2}}
\sup_{0\neq u\in E}\frac{h_n[u]}{\|u\|_2^2}\\
&\geq\inf_{\substack{E\subset D(h)\\ \dim E=2}}\sup_{0\neq u\in E}\frac{h[u]}{\|u\|_2^2}=\lambda_2^{\rm rad}(H)=0.
\end{aligned}
\end{equation*}
Moreover, since $C_c^\infty(\R^d\setminus\{0\})_{\rm rad}$ is dense in $D(h)$, there exist two sequences $\{v_k\},\{z_k\}\subset C_c^\infty(\R^d\setminus\{0\})_{\rm rad}$ such that
\begin{equation*}
v_k\to v\text{ \ and \ }
z_k\to z\text{ \ in \ }D(h)\text{ \ as \ }k\to\infty.
\end{equation*}
Define $E_k:=\operatorname{span}\{v_k,z_k\}$ for all sufficiently large $k>0$.
We then have
\begin{equation*}
\lim_{k\to\infty}
\sup_{\substack{u\in E_k\\ \|u\|_2=1}}h[u]=\sup_{\substack{u\in E\\ \|u\|_2=1}}h[u]=0.
\end{equation*}
For every fixed $k\in\N^+$, $W$ is bounded on $\operatorname{supp}v_k\cup\operatorname{supp}z_k$. This then implies that for all sufficiently large $n\in\N^+$,
\begin{equation*}
W_n\equiv W
\ \ \text{on }
\operatorname{supp}v_k\cup\operatorname{supp}z_k.
\end{equation*}
Consequently, we have $h_n[u]=h[u]$ for $u\in E_k$. The min--max principle therefore implies that
\begin{equation*}
\limsup_{n\to\infty}\lambda_{2,n}\leq
\sup_{\substack{u\in E_k\\ \|u\|_2=1}}h[u].
\end{equation*}
Setting $k\to\infty$ and using $\lambda_{2,n}\geq0$, we obtain that
\begin{equation}\label{eq:secondcon}
\lambda_{2,n}\to0\text{ \ as \ }n\to\infty.
\end{equation}

Since $h[v]<0$, we can fix  sufficiently large $k\in\N^+$ such that $h[v_k]<0$. For all sufficiently large $n\in\N^+$, we then have $h_n[v_k]=h[v_k]<0$, and hence
$\lambda_{1,n}\leq\frac{h_n[v_k]}{\|v_k\|_2^2}<0\leq\lambda_{2,n}$.
Note that $W_n\in L^\infty(\R^d)$ and
$W_n(x)\to0$ as $|x|\to\infty$. We thus obtain that $W_n\bigl((-\Delta)^s+1\bigr)^{-1}$ is compact on $L^2(\R^d)$
by the local Rellich compactness. Hence, by
Weyl's theorem \cite[Theorem~XIII.14]{ReedSimonIV}, we have
\begin{equation*}
\sigma_{\rm ess}(H_n)=\sigma_{\rm ess}\Big((-\Delta)^s+1\Big)=[1,\infty).
\end{equation*}
It thus follows from \eqref{eq:secondcon} that $\lambda_{1,n}
<\lambda_{2,n}<1=\inf\sigma_{\rm ess}(H_n)$ holds
for all sufficiently large $n\in\N^+$.
Set $V_n:=1-W_n$. It then follows from \eqref{eq:smoothcut} and Lemma \ref{lem3.1} that
\begin{equation*}
V_n\in C^{0,1}(\R^d)\cap L^\infty(\R^d),\ \  V_n'(r)\geq0.
\end{equation*}
Consequently, by \cite[Theorem~2.3]{FLS2016}, there exist a normalized second radial eigenfunction $Z_n$ satisfying $\|Z_n\|_2=1$ and $r_n\in(0,\infty)$ such that, after replacing $Z_n$ by $-Z_n$ if necessary,
\begin{equation}\label{eq:onenode}
H_nZ_n=\lambda_{2,n}Z_n\text{ \ in \ }\R^d, \ \ Z_n(r)\geq0\ (0<r<r_n),\ \ Z_n(r)\leq0\ (r>r_n).
\end{equation}
Moreover, since $h\leq h_n$, we have
\begin{equation*}
\mathfrak q_{s,\tau}[Z_n]+1-(p-1)\int_{\R^d}Q^{p-2}|Z_n|^2 dx\leq h_n[Z_n] = \lambda_{2,n}.
\end{equation*}
Using \eqref{eq:nonforms} and \eqref{eq:secondcon}, where $0<\varepsilon<(p-1)^{-1}$, we obtain that
\begin{equation}\label{eq:unifbound}
\sup_n\mathfrak q_{s,\tau}[Z_n]<\infty.
\end{equation}
It then follows from \eqref{2.3} that for every $0<\sigma<s$,
\begin{equation}
\sup_n\|Z_n\|_{H^\sigma}<\infty.
\label{eq:Hsigma}
\end{equation}

Let $\eta_R\in C^\infty(\R^d)$ satisfy $0\leq\eta_R\leq1$ and
\begin{equation*}
  \eta_R=0\text{ \ on \ }B_R(0),\ \ \eta_R=1\text{ \ on \ }\R^d\setminus B_{2R}(0),\ \ |\nabla\eta_R|\leq C/R.
\end{equation*}
It then yields that for sufficiently large $R>0$,
\begin{equation}\label{eq:IMS}
\left|\langle(-\Delta)^sZ_n,\eta_R^2Z_n\rangle-\|(-\Delta)^{s/2}(\eta_RZ_n)\|_2^2\right|
\leq CR^{-2s}\|Z_n\|_2^2.
\end{equation}
Testing \eqref{eq:onenode} with $\eta_R^2Z_n$ and using $W_n\leq W$, we derive from \eqref{eq:IMS} that
\begin{equation}\label{eq:IMStail}
\big(1-\lambda_{2,n}-\sup_{r\geq R}W(r)\big)\|\eta_RZ_n\|_2^2\leq CR^{-2s}.
\end{equation}
Since $\lim_{n\to\infty}\lambda_{2,n}=0$ and $\lim_{R\to\infty}\sup_{r\geq R}W(r)=0$,
it follows from \eqref{eq:IMStail} that
\begin{equation*}
\sup_{n\geq N_R}\int_{|x|\geq2R}|Z_n|^2 dx\leq CR^{-2s}.
\end{equation*}
By the local compactness of the embedding
$H^\sigma(\R^d)\hookrightarrow L^2_{\rm loc}(\R^d)$ for any fixed $\sigma\in(0,s)$,
and the uniform tail estimate above,  we obtain from \eqref{eq:Hsigma} that there exist $Z_*\in L^2(\R^d)$ and a subsequence, still denoted by $\{Z_n\}$, of $\{Z_n\}$ such that
\begin{equation}\label{eq:strongL2}
Z_n\to Z_*\text{ \ strongly in \ }L^2(\R^d)\text{ \ as  \ }n\to\infty,
\end{equation}
which implies that $ \|Z_*\|_2=1$.

It follows from \eqref{eq:unifbound} that
$\sup_n\|Z_n\|_{D(h)}<\infty$. This then implies that up to a subsequence if necessary,
\begin{equation*}
Z_n\rightharpoonup Z_*\text{ \ in \ }D(h)\text{ \ as \ }n\to\infty.
\end{equation*}
For any $\phi\in C_c^\infty(\R^d\setminus\{0\})_{\rm rad}$, since $W$ is bounded on $\operatorname{supp}\phi$, we have $W_n\equiv W$ on $\operatorname{supp}\phi$ for all sufficiently large $n\in\N^+$. This thus implies that
\begin{equation*}
h[Z_n,\phi]=h_n[Z_n,\phi]=\lambda_{2,n}\langle Z_n,\phi\rangle\to0\text{ \ as \ }n\to\infty.
\end{equation*}
Since $Z_n\rightharpoonup Z_*$ weakly in $D(h)$ as $n\to\infty$, we also have $h[Z_n,\phi]\to h[Z_*,\phi]$ as $n\to\infty$. Thus, $h[Z_*,\phi]=0$.
Since
$C_c^\infty(\R^d\setminus\{0\})_{\rm rad}$ is a form core of $h$, this
identity holds for every $\phi\in D(h)$, and hence $Z_*\in\ker H$.

Up to a subsequence if necessary, we now have
\begin{equation*}
r_n\to r_*\in[0,\infty]\text{ \ as \ }n\to\infty,
\end{equation*}
where $r_n\in\R^+$ is as in \eqref{eq:onenode}.  If $r_*=0$, then for every $0<a<b<\infty$,
\begin{equation*}
Z_n\leq0\text{ \ on \ }\{a<|x|<b\}\text{ \ for all sufficiently large \ }n\in\N^+.
\end{equation*}
 Hence
\begin{equation*}
\|(Z_*)_+\|_{L^2(\{a<|x|<b\})}\leq\|Z_n-Z_*\|_2\to0\text{ \ as \ }n\to\infty.
\end{equation*}
Setting $a\downarrow0$ and $b\uparrow\infty$, it gives that $Z_*\leq0~a.e.$ in $\R^d$, which then contradicts \eqref{eq:1moment} and \eqref{eq:strongL2}. Similarly, if $r_*=\infty$, then it implies that $Z_*\geq0~a.e.$ in $\R^d$.
This also contradicts \eqref{eq:1moment} and \eqref{eq:strongL2}. We thus conclude that $0<r_*<\infty$.

It follows from \eqref{eq:onenode} and \eqref{eq:strongL2} that
\begin{equation}\label{eq:lim1node}
Z_*(r)\geq0\quad\text{for }~a.e.~0<r<r_* ,\ \  Z_*(r)\leq0\quad\text{for }~a.e.~ r>r_*.
\end{equation}
Moreover, it follows from \eqref{eq:1moment} that
$Z_*$ is nonzero on the sets around $r_*$ with positive measure. Since $Q>0$ is strictly decreasing in $r>0$, it follows from
\eqref{eq:lim1node} that
\begin{equation*}
Q(r)\Big[Q(r)^{p-2}-Q(r_*)^{p-2}\Big]Z_*(r)\geq0\text{ \ for   \ } a.e. \ \ r>0,
\end{equation*}
where the strict inequality holds on a set of positive measure. We therefore deduce that
\begin{align*}
0&<\int_{\R^d}Q\bigl(Q^{p-2}-Q(r_*)^{p-2}\bigr)Z_* dx\\
&=\int_{\R^d}Q^{p-1}Z_* dx-Q(r_*)^{p-2}\int_{\R^d}QZ_* dx=0,
\end{align*}
a contradiction,
where the last identity follows from
\eqref{eq:1moment} and
\eqref{eq:2moment}. This proves that
$\ker H\big|_{\rm rad}=\{0\}$. We thus complete the proof of Proposition \ref{prop:radialnon}.
\qed

\subsection{Full-space nondegeneracy  }\label{rigang:sec-full-space}

Recall that
Proposition~\ref{prop:radialnon}  excludes
the radial zero modes of $L_{+,\tau}$. To complete the proof of
\eqref{1.14}, in this subsection we further exclude the nonradial zero
modes of the linearized operator.

Throughout this section, we set
\begin{equation*}
d\geq3,\ \  0<s<1,\ \  \alpha:=2s,\ \ 2<p<2_s^*:=\frac{2d}{d-2s}.
\end{equation*}
For each $0\leq\tau\leq1$, let $Q(x):=Q_\tau(x)=Q_\tau(|x|)>0$ be a positive radial weak solution of
\begin{equation}\label{eq:Qtau}
\bigl(( -\Delta)^s-\tau\Cc_{d,s}|x|^{-\alpha}\bigr)Q+Q=Q^{p-1}\text{ \ in \ }\R^d.
\end{equation}
Recall from \eqref{2.2} that
\begin{equation}\label{eq:htau}
\mathfrak h_{s,\tau}[u]=\mathfrak q_{s,\tau}[u]+\|u\|_2^2,
\end{equation}
and
\begin{equation*}
D(\mathfrak h_{s,\tau})=
\begin{cases}
H^s(\R^d),&0\leq\tau<1;\\
\mathcal Q^s,&\tau=1.
\end{cases}
\end{equation*}
Moreover, we obtain from \eqref{eq:nonforms} that
for every fixed $\varepsilon>0$,
there exists a constant $C_\varepsilon>0$ such that
$$
\int_{\R^d}Q^{p-2}|u|^2\,dx
\leq \varepsilon\mathfrak q_{s,\tau}[u]
+C_\varepsilon\|u\|_2^2.$$
Thus, the KLMN theorem  (cf. \cite[Theorem 3.19]{Lewin}) gives the closed and lower-semibounded form
\begin{equation}\label{eq:linearizedform}
\mathfrak l_{+,\tau}[u]:=\mathfrak h_{s,\tau}[u]-(p-1)\int_{\R^d}Q^{p-2}|u|^2 dx,
\ \ D(\mathfrak l_{+,\tau})=D(\mathfrak h_{s,\tau}),
\end{equation}
for which we define
\begin{equation}\label{2.24A}
\mathfrak l_{+,\tau}(u,v)
:=\mathfrak h_{s,\tau}[u,v]
-(p-1)\int_{\R^d}Q^{p-2}uv\,dx,
\quad u,v\in D(\mathfrak l_{+,\tau}).
\end{equation}
We denote the associated self-adjoint operator by
\begin{equation}\label{eq:Lplus}
L_{+,\tau}:=(-\Delta)^s-\tau\Cc_{d,s}|x|^{-\alpha}+1-(p-1)Q^{p-2}\text{ \ in \ }\R^d.
\end{equation}

Let $\{Y_{\ell,m}\}$ be a real orthonormal basis of spherical harmonics.
Define for $r>0$ and $\omega\in\mathbb S^{d-1}$,
\begin{equation*}
(P_{\ell,m}u)(r\omega):=Y_{\ell,m}(\omega)\int_{\mathbb S^{d-1}}u(r\nu)Y_{\ell,m}(\nu) d\nu,
\ \ V_{\ell,m}:=\operatorname{Ran}P_{\ell,m}.
\end{equation*}
These are the orthogonal angular projections in $L^2(\R^d)$. For each $\ell\geq0$, we define the angular momentum channel
of degree $\ell$ by
\begin{equation*}
V_\ell:=\bigoplus_m V_{\ell,m},
\end{equation*}
where the sum runs over all spherical harmonics of degree $\ell\geq0$.
The next lemma extends the angular decomposition to the critical
closed form domain, without separating the kinetic and Hardy energies.

\begin{Lemma}
\label{lem:directsum}
Each angular projection $P_{\ell,m}$ maps $\mathcal Q^s$ into itself and satisfies
\begin{equation}\label{eq:formcomm}
\mathfrak l_{+,1}(P_{\ell,m}u,v)=\mathfrak l_{+,1}(u,P_{\ell,m}v),
\ \ \ \forall \  u,v\in\mathcal Q^s,
\end{equation}
where $l_{+,1}(\cdot,\cdot)$ is defined by \eqref{2.24A}.
In particular, we have for $(\ell,m)\neq(\ell',m')$,
\begin{equation}\label{eq:5.02}
\mathfrak l_{+,1}(P_{\ell,m}u,P_{\ell',m'}v)=0,
\ \ \ \forall \  u,v\in\mathcal Q^s.
\end{equation}
\end{Lemma}

\noindent\textbf{Proof.}
Fix $P=P_{\ell,m}$. By \cite[Lemma~1.2]{BogdanMerz2024}, we have
\begin{equation}\label{eq:5.01}
\mathfrak h_{s,1}(Pu,v)=\mathfrak h_{s,1}(u,Pv),
\ \ \ \forall \   u,v\in C_c^\infty(\mathbb R^d\setminus\{0\}).
\end{equation}
Since $P$ preserves $ C_c^\infty(\mathbb R^d\setminus\{0\})$, it follows from \eqref{eq:5.01} that
\begin{equation*}
\mathfrak h_{s,1}[u]
=
\mathfrak h_{s,1}[Pu]+\mathfrak h_{s,1}[(I-P)u]
\geq \mathfrak h_{s,1}[Pu],
\ \ \ \forall \   u\in C_c^\infty(\mathbb R^d\setminus\{0\}).
\end{equation*}
By the density of $ C_c^\infty(\mathbb R^d\setminus\{0\})$ in $\mathcal Q^s$, this implies that $P$ maps $\mathcal Q^s$ into itself.
Since $Q>0$ is radial, we have
\begin{equation*}
\int_{\R^d}Q^{p-2}(Pu) v dx=\int_{\R^d}Q^{p-2}u{Pv} dx,\ \ \ \forall \   u,v\in C_c^\infty(\mathbb R^d\setminus\{0\}).
\end{equation*}
Hence, \eqref{eq:formcomm} holds in $ C_c^\infty(\mathbb R^d\setminus\{0\})$. By \eqref{eq:nonforms}, $\mathfrak l_{+,1}$ is continuous on $\mathcal Q^s\times\mathcal Q^s$. Using the density of $ C_c^\infty(\mathbb R^d\setminus\{0\})$ in $\mathcal {Q}^s$, we can then derive that \eqref{eq:formcomm} holds for all $u,v\in\mathcal Q^s$. Taking $v=P_{\ell',m'}v$ into \eqref{eq:formcomm}, it thus yields \eqref{eq:5.02}. We thus complete the proof of Lemma~\ref{lem:directsum}.
\qed

\vspace{0.5em}

For $0<\tau<1$, since $D(\mathfrak l_{+,\tau})=H^s(\R^d)$, the same angular reduction follows directly from the
rotational invariance of $(-\Delta)^s$ and the radiality of
$|x|^{-2s}$ and $Q^{p-2}$, while the critical case $\tau=1$ follows
from Lemma~\ref{lem:directsum}. Hence, each
$V_{\ell,m}$ reduces to $L_{+,\tau}$ for every $0<\tau\leq1$. Let
$L_{+,\tau}^{(\ell,m)}$ be the corresponding self-adjoint restriction
of $L_{+,\tau}$ on $V_{\ell,m}$. Then
\begin{equation*}
L_{+,\tau}=\bigoplus_{\ell,m}L_{+,\tau}^{(\ell,m)},
\quad
\ker L_{+,\tau}
=\bigoplus_{\ell,m}\ker L_{+,\tau}^{(\ell,m)}.
\end{equation*}
Thus, to reduce the full-space nondegeneracy of $L_{+,\tau}$ to the radial nondegeneracy of $L_{+,\tau}$, it next suffices to exclude zero modes of $L_{+,\tau}$ in all angular channels with $\ell\geq1$.

We now consider the first momentum angular channel of $L_{+,\tau}$. Take $Y_{1,m}(\omega)=a_m\cdot\omega$ satisfying
$\|Y_{1,m}\|_{L^2(\mathbb S^{d-1})}=1$, and use the solid-harmonic
coordinate
\begin{equation}\label{eq:solidcoor}
[f]_{1,m}(x):=f(r)rY_{1,m}(\omega)=f(r)a_m\cdot x,\ \  r=|x|,\quad\omega=x/r,
\end{equation}
where the corresponding radial measure is
\begin{equation}\label{eq:solidmea}
dm_1(r):=r^{d+1} dr,\ \ \|[f]_{1,m}\|_2^2=\int_0^\infty|f|^2 dm_1.
\end{equation}
Taking $\ell=1$ and $\sigma=-1$, it then yields from  \cite[Theorem~1.5]{BogdanMerz2024} that for any
$f\in C_c^\infty(0,\infty)$,
\begin{equation}\label{eq:1channel}
\|(-\Delta)^{s/2}[f]_{1,m}\|_2^2=\frac12\iint_{(0,\infty)^2}|f(r)-f(\rho)|^2K_1(r,\rho) dm_1(r)dm_1(\rho):=\mathcal E_1(f),
\end{equation}
where $K_1(r,\rho):=\nu_{(d+1)/2}(r,\rho)$ is strictly positive, and satisfies the estimate
\begin{equation}\label{eq:K1estimate}
K_1(r,\rho)\sim
\frac{1}{(r+\rho)^{d+1}|r-\rho|^{1+2s}},
\ \  r\neq\rho,
\end{equation}
in view of \cite[Proposition~2.4]{BogdanMerz2024}.
 Define for $f,g\in C_c^\infty(0,\infty)$,
\begin{equation}\label{eq:5.12}
\mathcal E_1(f,g):=\frac12\int_0^\infty \int_0^\infty\big(f(r)-f(\rho)\big)\big(g(r)-g(\rho)\big)K_1(r,\rho) dm_1(r)dm_1(\rho).
\end{equation}
In the following, we derive the weak equation of $-Q'/r$ in the first angular momentum channel.
\begin{Lemma} \label{rigang:lem-first-channel-differentiation}
Suppose $Q(x):=Q_\tau>0$ is a positive radial ground state of \eqref{eq:Qtau} with $0\leq\tau\leq1$.
Set
\begin{equation*}
\psi(r):=-\frac{Q'(r)}r\geq0,\ \ g_\tau(r):=\alpha\tau\Cc_{d,s}r^{-\alpha-2}Q(r),\ \  r>0,
\end{equation*}
and
\begin{equation*}
V_\tau(r):=1-\tau\Cc_{d,s}r^{-\alpha}-(p-1)Q(r)^{p-2},\ \  r>0.
\end{equation*}
Then we have
\begin{equation}\label{eq:weaksup}
\mathcal E_1(\psi,\eta)+\int_0^\infty V_\tau\psi\eta dm_1=\int_0^\infty g_\tau\eta dm_1,\ \ \ \forall \ \eta\in C_c^\infty(0,\infty),
\end{equation}
where the integral $\E_1(\psi,\eta)$ is absolutely convergent.
\end{Lemma}

\noindent\textbf{Proof.}
Define $U_m:=[\psi]_{1,m}=-a_m\cdot\nabla Q$ in $\R^d\setminus\{0\}$. Since $Q\in L^p(\R^d)$ is nonincreasing in $|x|>0$, we have
\begin{equation*}
Q(r)\leq C\|Q\|_p r^{-d/p},\quad r>0.
\end{equation*}
It then yields that
\begin{equation*}
\int_0^1(-Q'(r))r^{d-1} dr<\infty,\quad\int_1^\infty(-Q'(r))r^{-1-2s} dr\leq Q(1).
\end{equation*}
We thus obtain that
\begin{equation}\label{eq:Ulocal1}
\int_{|x|<1}|U_m(x)| dx\leq C\int_0^1(-Q'(r))r^{d-1} dr<\infty,
\end{equation}
and
\begin{equation}\label{eq:weightedL1}
\int_{|x|>1}\frac{|U_m(x)|}{|x|^{d+2s}} dx\leq C\int_1^\infty(-Q'(r))r^{-1-2s} dr<\infty.
\end{equation}
Since $r^{d-1}Q(r)\to0$ as $r\searrow0$, we obtain that $U_m=-a_m\cdot\nabla Q$ in $D'(\R^d\setminus\{0\})$. Differentiating \eqref{eq:Qtau}, we then derive that
\begin{equation}\label{eq:dissup}
L_{+,\tau}U_m=[g_\tau]_{1,m}\ \ \ \ \text{in}\ \ \mathcal D'(\R^d\setminus\{0\}),
\end{equation}
where $g_\tau(r)=\alpha\tau\Cc_{d,s}r^{-\alpha-2}Q(r)$.

Fix $\eta\in C_c^\infty(0,\infty)$.
It yields that
\begin{equation}\label{A1}
\begin{aligned}
\langle(-\Delta)^sU_m,[\eta]_{1,m}\rangle
&=\int_{\R^d}U_m(x)(-\Delta)^s[\eta]_{1,m}(x) dx\\
&=\frac{c_{d,s}}2\iint_{\R^d\times\R^d}\frac{(U_m(x)-U_m(y))([\eta]_{1,m}(x)-[\eta]_{1,m}(y))}
{|x-y|^{d+\alpha}} dx dy,
\end{aligned}
\end{equation}
where $c_{d,s}>0$ is the normalization constant of the fractional
Laplacian. Since $\alpha<2$, using the local
smoothness of $U_m$ and $[\eta]_{1,m}$, we deduce from \eqref{eq:Ulocal1}
and \eqref{eq:weightedL1} that the integral \eqref{A1} is absolutely convergent.
Moreover, for $0<2\varepsilon<R$, choose $\chi_{\varepsilon,R}\in C_c^\infty\big((0,\infty),[0,1]\big)$ such that
\begin{equation}\label{eq:cutoff}
\chi_{\varepsilon,R}(r)\equiv1\ \ \text{for}\ \ 2\varepsilon\leq r\leq R,\ \ \chi_{\varepsilon,R}(r)\equiv0\ \ \text{for}\ \ 0<r\leq\varepsilon\ \ \text{and}\ \ r\geq2R.
\end{equation}
It then implies from \eqref{eq:5.12} that
\begin{equation*}
\begin{aligned}
\mathcal E_1(\chi_{\varepsilon,R}\psi,\eta)&=\left\langle(-\Delta)^s[\chi_{\varepsilon,R}\psi]_{1,m},[\eta]_{1,m}\right\rangle\\
&=\int_{\R^d}\chi_{\varepsilon,R}(|x|)U_m(x)(-\Delta)^s[\eta]_{1,m}(x) dx.
\end{aligned}
\end{equation*}
Since $[\eta]_{1,m}\in C_c^\infty(\R^d\setminus\{0\})$,
we have
\begin{equation*}
\bigl|(-\Delta)^s[\eta]_{1,m}(x)\bigr|\leq C_\eta(1+|x|)^{-d-\alpha}\ \ \ \ \text{in}\ \ \R^d.
\end{equation*}
By dominated convergence theorem, we thus obtain from
\eqref{eq:Ulocal1} and
\eqref{eq:weightedL1} that
\begin{equation}\label{eq:5.165}
\lim_{\substack{\varepsilon\searrow0\\R\nearrow\infty}}\mathcal E_1(\chi_{\varepsilon,R}\psi,\eta)=\langle(-\Delta)^sU_m,[\eta]_{1,m}\rangle.
\end{equation}

On the other hand, it yields that
\begin{equation}\label{eq:smallpaitail}
\int_0^{2\varepsilon}\psi(\rho)\rho^{d+1} d\rho=\int_0^{2\varepsilon}(-Q'(\rho))\rho^d d\rho
=O(\varepsilon^{d-d/p})\ \ \ \ \text{as}\ \ \varepsilon\searrow0,
\end{equation}
and
\begin{equation}\label{eq:largepsitail}
\int_R^\infty\psi(\rho)\rho^{-1-\alpha} d\rho
=\int_R^\infty(-Q'(\rho))\rho^{-2-\alpha} d\rho
\leq Q(R)R^{-2-\alpha}.
\end{equation}
Together with the local smoothness and $\alpha<2$, these estimates imply from \eqref{eq:K1estimate}
that $\mathcal E_1(\psi,\eta)$ is absolutely convergent. For sufficiently small $\varepsilon>0$ and large $R>0$, $\chi_{\varepsilon,R}\equiv1$ holds on a neighborhood of $\operatorname{supp}\eta$. We thus obtain from \eqref{eq:K1estimate} that for sufficiently large $R>0$,
\begin{equation*}
\begin{aligned}
&\bigl|\mathcal E_1(\chi_{\varepsilon,R}\psi,\eta)-\mathcal E_1(\psi,\eta)\bigr|=\bigl|\mathcal E_1((\chi_{\varepsilon,R}-1)\psi,\eta)\bigr|\\
\leq& C_\eta\left(\int_0^{2\varepsilon}\psi(\rho)\rho^{d+1} d\rho+\int_R^\infty\psi(\rho)\rho^{-1-\alpha} d\rho\right)\to0\ \ \text{as}\ \ \varepsilon\searrow0.
\end{aligned}
\end{equation*}
Together with \eqref{eq:5.165}, this implies that
\begin{equation}\label{eq:5.16}
\mathcal E_1(\psi,\eta)=\langle(-\Delta)^sU_m,[\eta]_{1,m}\rangle.
\end{equation}
Multiplying \eqref{eq:dissup} by $[\eta]_{1,m}$ and integrating over $\R^d$, we conclude from \eqref{eq:5.16} that \eqref{eq:weaksup} holds true.
We thus complete the proof of Lemma~\ref{rigang:lem-first-channel-differentiation}.
\qed

 \vspace{0.5em}

Applying Lemma \ref{rigang:lem-first-channel-differentiation}, we then prove the strict negativity of $Q'$.
\begin{Lemma} \label{lem:strict}
Suppose $Q(x):=Q_\tau$ is a positive radial ground state of \eqref{eq:Qtau} with $0\leq\tau\leq1$. Then it holds that for every $0\leq\tau\leq1$,
\begin{equation*}
\psi(r)=-\frac{Q'(r)}r>0,\   \ r>0.
\end{equation*}
\end{Lemma}

\noindent\textbf{Proof.}
Fix $0\leq\tau\leq1$ and suppose that $\psi(r_0)=0$ holds for some $r_0>0$. Since $\psi\geq0$ and $\psi\in C^\infty(0,\infty)$, we have
\begin{equation*}
\psi'(r_0)=0,\ \ \psi(\rho)\leq C_{r_0}|\rho-r_0|^2\ \ \text{for}\ \ |\rho-r_0|<r_0/2.
\end{equation*}
We thus obtain from \eqref{eq:K1estimate}, \eqref{eq:smallpaitail} and \eqref{eq:largepsitail} that
\begin{equation}\label{eq:T1absolute}
\begin{aligned}
&\int_0^\infty
|\psi(r_0)-\psi(\rho)|K_1(r_0,\rho) dm_1(\rho)\\
\leq& C_{r_0}\left(\int_0^{r_0/2}\psi(\rho)\rho^{d+1} d\rho+\int_{r_0/2}^{3r_0/2}|\rho-r_0|^{1-\alpha} d\rho+\int_{3r_0/2}^\infty\psi(\rho)\rho^{-1-\alpha} d\rho\right)<\infty.
\end{aligned}
\end{equation}
Set $U_m:=[\psi]_{1,m}$ and choose $x_0=r_0e$, where $e\in\mathbb S^{d-1}$ and $a_m\cdot e\neq0$. It then follows that $U_m(x_0)=0$ and $\nabla U_m(x_0)=0$.
For $\chi_{\varepsilon,R}(r)$ as in \eqref{eq:cutoff}, denote
\begin{equation*}
U_{m,\varepsilon,R}:=[\chi_{\varepsilon,R}\psi]_{1,m}\in C_c^\infty(\R^d\setminus\{0\}).
\end{equation*}
Fix a neighborhood $N\Subset\R^d\setminus\{0\}$ of $x_0$. For sufficiently small $\varepsilon>0$ and large $R>0$, we have $U_{m,\varepsilon,R}\equiv U_m$ in a neighborhood of $\overline N$.
Furthermore, we derive from \eqref{eq:Ulocal1} and \eqref{eq:weightedL1} that for sufficiently large $R>0$,
\begin{equation}\label{eq:pointwiseerror}
\begin{aligned}
&\sup_{x\in N}
\left|(-\Delta)^s(U_m-U_{m,\varepsilon,R})(x)\right|\\
=&c_{d,s}\sup_{x\in N}\left|\int_{\R^d}\frac{(1-\chi_{\varepsilon,R}(|y|))U_m(y)}
{|x-y|^{d+\alpha}} dy\right|\\
\leq&C_N\left(\int_0^{2\varepsilon}\psi(\rho)\rho^d d\rho+\int_R^\infty\psi(\rho)\rho^{-\alpha} d\rho\right)\\
\to&0\ \ \text{as}\ \ \varepsilon\searrow0.
\end{aligned}
\end{equation}
Since $(-\Delta)^sU_{m,\varepsilon,R}$ is continuous on $N$, the uniform convergence of \eqref{eq:pointwiseerror} implies that $(-\Delta)^sU_m$ is also continuous on $N$. This then implies that \eqref{eq:dissup} holds pointwise on $N$.

It follows from \eqref{eq:solidcoor} that
\begin{equation*}
U_{m,\varepsilon,R}(x)
=[\chi_{\varepsilon,R}\psi]_{1,m}(x)
=\chi_{\varepsilon,R}(|x|)\psi(|x|)(a_m\cdot x)
\  \ \ \ \text{in}\ \ \R^d\setminus\{0\}.
\end{equation*}
Since $\psi(r_0)=\psi'(r_0)=0$ and $\alpha<2$, $(-\Delta)^sU_{m,\varepsilon,R}(x_0)$ is well-defined. Using the polar coordinates, we then obtain from \cite[Proposition~4.3(3)]{BogdanMerz2024} with $\ell=1$ that
\begin{equation*}
\frac{(-\Delta)^sU_{m,\varepsilon,R}(x_0)}{a_m\cdot x_0}=-\int_0^\infty\chi_{\varepsilon,R}(\rho)\psi(\rho)K_1(r_0,\rho) dm_1(\rho).
\end{equation*}
By dominated convergence theorem, we deduce from  \eqref{eq:T1absolute} that for sufficiently large $R>0$,
\begin{equation*}
-\int_0^\infty\chi_{\varepsilon,R}(\rho)\psi(\rho)K_1(r_0,\rho) dm_1(\rho)\to\mathcal T_1\psi(r_0)\ \ \text{as}\ \ \varepsilon\searrow0,
\end{equation*}
where
$\mathcal T_1\psi(r_0)
:=\int_0^\infty
(\psi(r_0)-\psi(\rho))K_1(r_0,\rho) dm_1(\rho)$.
Together with \eqref{eq:pointwiseerror},
this implies that
\begin{equation*}
(-\Delta)^sU_m(x_0)
=(a_m\cdot x_0)\mathcal T_1\psi(r_0).
\end{equation*}
Taking $x=x_0$ in \eqref{eq:dissup} and dividing by $a_m\cdot x_0$, we then obtain that
\begin{equation}\label{eq:point1channel}
\mathcal T_1\psi(r_0)+V_\tau(r_0)\psi(r_0)
=g_\tau(r_0)\geq0.
\end{equation}

Since $Q(|x|)$ is strictly decreasing in $|x|>0$, we have $0\leq\psi\not\equiv0$ on $(0,\infty)$.
By the continuity of $\psi$, there exists a nonempty open
interval $I\subset(0,\infty)$ such that $\psi>0$ on $I$.
We then obtain that
\begin{equation*}
\begin{aligned}
\mathcal T_1\psi(r_0)
&=-\int_0^\infty
\psi(\rho)K_1(r_0,\rho) dm_1(\rho)
 \leq-\int_I
\psi(\rho)K_1(r_0,\rho) dm_1(\rho)<0,
\end{aligned}
\end{equation*}
where $K_1>0$ is as in \eqref{eq:1channel}.
This however contradicts \eqref{eq:point1channel} in view of $\psi(r_0)=0$. We therefore conclude that $\psi(r)>0$ holds for every $r>0$ and every $0\leq\tau\leq1$. This thus completes the proof of Lemma~\ref{lem:strict}.
\qed

For every $f\in C_c^\infty(0,\infty)$, it implies from Lemma~\ref{lem:strict} that $f^2/\psi\in C_c^\infty(0,\infty)$. Using the identity
\begin{equation*}
\begin{aligned}
&(f(r)-f(\rho))^2-(\psi(r)-\psi(\rho))\left(\frac{f(r)^2}{\psi(r)}-\frac{f(\rho)^2}{\psi(\rho)}\right)\\
=& \psi(r)\psi(\rho)\left(\frac{f(r)}{\psi(r)}-\frac{f(\rho)}{\psi(\rho)}\right)^2,\quad \forall~ r,\rho>0,
\end{aligned}
\end{equation*}
and taking $\eta=f^2/\psi$ in \eqref{eq:weaksup}, we then derive from \eqref{eq:linearizedform} and \eqref{eq:1channel} that
\begin{equation}\label{eq:1cinden}
\begin{aligned}
 \mathfrak l_{+,\tau}[[f]_{1,m}]
=&\frac12\iint_{(0,\infty)^2}K_1(r,\rho)\psi(r)\psi(\rho)\left(\frac{f(r)}{\psi(r)}-\frac{f(\rho)}{\psi(\rho)}\right)^2
 dm_1(r)dm_1(\rho)\\
&+\int_0^\infty\frac{g_\tau(r)}{\psi(r)}f(r)^2 dm_1(r).
\end{aligned}
\end{equation}
The following Proposition \ref{prop:1channel} describes the kernel of $L_{+,\tau}$ on $V_1$.

\begin{Proposition} \label{prop:1channel}
For every $0\leq\tau\leq1$, the operator $L_{+,\tau}$ is nonnegative on $V_1$, where $L_{+,\tau}$ is as in \eqref{eq:Lplus}. Moreover, we have
\begin{equation*}
\ker\bigl(L_{+,\tau}|_{V_1}\bigr)=
\begin{cases}
\{0\},&\text{if}\ \ 0<\tau\leq1;\\
\operatorname{span}\{\partial_jQ:~j=1,\cdots,d\},&\text{if}\ \ \tau=0.
\end{cases}
\end{equation*}
\end{Proposition}

\noindent\textbf{Proof.}
Since the set $\{[f]_{1,m}:f\in C_c^\infty(0,\infty)\}$ is dense in $D(\mathfrak l_{+,\tau})\cap V_{1,m}$ with respect to the norm $\mathfrak h_{s,\tau}^{1/2}$, we obtain from \eqref{eq:1cinden} that
\begin{equation*}
\mathfrak l_{+,\tau}[F]\geq0,\ \  F\in D(\mathfrak l_{+,\tau})\cap V_{1,m}.
\end{equation*}
Let $F=[f]_{1,m}\in D(\mathfrak l_{+,\tau})\cap V_{1,m}$ satisfy that $\mathfrak l_{+,\tau}[F]=0$. Then there exists a sequence $\{f_n\}\subset C_c^\infty(0,\infty)$ such that $[f_n]_{1,m}\to F$ strongly in $D(\mathfrak l_{+,\tau})\cap V_{1,m}$ as $n\to\infty$. We thus deduce from \eqref{eq:solidmea} that
\begin{equation}\label{eq:5.27}
f_n\to f\ \ \text{strongly in}\ \ L^2((0,\infty),dm_1)\ \ \text{as}\ \ n\to\infty,\ \ \mathfrak l_{+,\tau}[[f_n]_{1,m}]\to0\ \ \text{as}\ \ n\to\infty.
\end{equation}
For $0<\tau\leq1$, the function $g_\tau/\psi$ is continuous and strictly positive on every compact interval $K\Subset(0,\infty)$. We thus get from \eqref{eq:1cinden} that
\begin{equation*}
\int_K f_n^2 dm_1\leq\left(\min_K\frac{g_\tau}{\psi}\right)^{-1}\mathfrak l_{+,\tau}[[f_n]_{1,m}]\to0\ \ \text{as}\ \ n\to\infty.
\end{equation*}
This then implies from \eqref{eq:5.27} that $f\equiv0$ holds on every such interval $K$. We therefore obtain that $F\equiv0$ in $L^2(\R^d)$. Since $m$ is arbitrary, the kernel of the restriction of $L_{+,\tau}$ to $V_1$ is $\{0\}$ for $0<\tau\leq1$.

The case $\tau=0$ follows directly from \cite{FLS2016}. We thus complete the proof of Proposition~\ref{prop:1channel}.
\qed

\vspace{0.5em}

We next compare the higher angular channels with the first channel of  $\mathfrak l_{+,\tau}$. Recall that $\mathfrak h_{s,\tau}$ and $\mathfrak l_{+,\tau}$ are defined in \eqref{2.2} and \eqref{eq:linearizedform}, respectively. Note from \cite[Proposition~1.4 and Section~1.4]{BogdanMerz2024} that for $\ell\geq1$, the sharp Hardy constant in $V_\ell$ is $\kappa_\ell>\Cc_{d,s}$. Taking the closures in each angular channel, we derive from \eqref{eq:nonforms} that
\begin{equation}\label{eq:nonform}
D(\mathfrak l_{+,\tau})\cap V_\ell=D(\mathfrak h_{s,\tau})\cap V_\ell=H^s(\R^d)\cap V_\ell,\ \  \ell\geq1.
\end{equation}
For each $\ell\geq1$, choose a spherical harmonic $Y_\ell$ of degree $\ell$ satisfying
$\|Y_\ell\|_{L^2(\mathbb S^{d-1})}=1$, and write
\begin{equation*}
u(x)=v(r)Y_\ell(\omega),\ \  r=|x|,\quad\omega=x/|x|.
\end{equation*}
Let $\mathfrak l_{+,\tau}^{(\ell)}$ denote the radial closed form on $\mathcal H_r:=L^2((0,\infty),r^{d-1} dr)$ given by
\begin{equation*}
\mathfrak l_{+,\tau}^{(\ell)}[v]:=\mathfrak l_{+,\tau}[u]=\|(-\Delta)^{s/2}u\|_2^2+\int_0^\infty V_\tau(r)|v(r)|^2r^{d-1} dr,
\end{equation*}
where $V_\tau(r):=1-\tau\Cc_{d,s}r^{-2s}-(p-1)Q_\tau(r)^{p-2}$,
and $Q_\tau>0$ is a ground state of \eqref{eq:Qtau}. By \eqref{eq:nonform}, the form domain of $\mathfrak l_{+,\tau}^{(\ell)}$ satisfies
\begin{equation*}
D(\mathfrak l_{+,\tau}^{(\ell)})=\bigl\{v\in\mathcal H_r:~~v(|x|)Y_\ell(x/|x|)\in H^s(\R^d)\bigr\}.
\end{equation*}

Applying \cite[Lemma~C.5 and Remark~C.6]{FLS2016} with zero potential and then adding the common potential $V_\tau$, we obtain that $D(\mathfrak l_{+,\tau}^{(\ell)}) \subset D(\mathfrak l_{+,\tau}^{(1)})$ holds for $\ell\geq2$. Together with Proposition~\ref{prop:1channel}, we obtain that
\begin{equation}\label{eq:fullcom}
\mathfrak l_{+,\tau}^{(\ell)}[v]>\mathfrak l_{+,\tau}^{(1)}[v]\geq0,\ \ 0\not\equiv v\in D(\mathfrak l_{+,\tau}^{(\ell)}),\quad\ell\geq2.
\end{equation}
It then follows from \eqref{eq:fullcom} that for every
$0\leq\tau\leq1$, the kernel of $L_{+,\tau}$
restricted to each angular channel $V_\ell$
with $\ell\geq2$ is also trivial.

\begin{Theorem} \label{thm:fullkernel}
Suppose $Q_\tau>0$ is a ground state of \eqref{eq:Qtau}, and let $L_{+,\tau}$ be defined by \eqref{eq:Lplus}, where $\tau\in[0,1]$. Then for $0<\tau\leq1$,
\begin{equation}\label{eq:fullkernel}
\ker L_{+,\tau}=\big\{0\big\}.
\end{equation}
At $\tau=0$, we have
\begin{equation}\label{eq:0coupling}
\ker L_{+,0}=\operatorname{span}\big\{\partial_1Q_0,\cdots,\partial_dQ_0\big\}.
\end{equation}
\end{Theorem}

\noindent\textbf{Proof.}
Applying Lemma~\ref{lem:directsum} with $\tau=1$ and the usual angular decomposition in $H^s(\R^d)$ for $\tau<1$, we deduce that
\begin{equation*}
\ker L_{+,\tau}=\bigoplus_{\ell\geq0}\ker\bigl(L_{+,\tau}|_{V_\ell}\bigr),\ \  0\leq\tau\leq1.
\end{equation*}
For $0<\tau\leq1$, Proposition~\ref{prop:radialnon} excludes the radial kernel of $L_{+,\tau}$, Proposition~\ref{prop:1channel} excludes the first-channel kernel of $L_{+,\tau}$, and \eqref{eq:fullcom} excludes the kernels in all higher angular channels of $L_{+,\tau}$.
Thus, we have
\begin{equation*}
\ker\bigl(L_{+,\tau}|_{V_\ell}\bigr)=\{0\},\ \  \ell\geq0,
\end{equation*}
which hence proves \eqref{eq:fullkernel}.

As for the case $\tau=0$, \eqref{eq:0coupling} follows directly from \cite[Theorem~3.3]{FLS2016}. This thus completes the proof of Theorem~\ref{thm:fullkernel}.
\qed

\section{Uniqueness of Ground States and Its Applications}\label{sec:cont-fixed-order}
In this section, we mainly prove the uniqueness of Theorem \ref{thm1.1} and address the applications of Theorem \ref{thm1.1}.
Throughout this section, we always suppose that
\begin{equation}\label{eq:admis}
d\geq3,\ \  0<s<1,\ \  2<p<2_s^*:=\frac{2d}{d-2s}.
\end{equation}
For $0\leq\tau\leq1$, set
\begin{equation}\label{cont:B-R-tau}
\mathcal{R}_\tau:=(A_{s,\tau}+1)^{-1}\ \ \text{on}\ \ L^2(\R^d),
\end{equation}
where the operator $A_{s,\tau}$ is as in \eqref{A}.
For $f\in L^{p'}_{\rm rad}(\R^d)$, where $p'=p/(p-1)$, we use the same notation $\mathcal{R}_\tau f$ for the unique radial function in $D(\mathfrak h_{s,\tau})$ satisfying
\begin{equation}\label{eq:varreso}
\mathfrak h_{s,\tau}[\mathcal{R}_\tau f,v]=\int_{\R^d}f v dx,\ \  v\in D(\mathfrak h_{s,\tau})\cap L^2_{\rm rad}(\R^d).
\end{equation}
The existence and uniqueness of $\mathcal{R}_\tau f$ follow from Lemma~\ref{lem2.1} and the Lax--Milgram theorem. By the uniqueness  of $\mathcal{R}_\tau f$, this definition agrees with $(A_{s,\tau}+1)^{-1}f$ for any $f\in L^2_{\rm rad}(\R^d)\cap L^{p'}_{\rm rad}(\R^d)$. In particular, a positive radial function $Q\in L^p(\R^d)$ is a weak solution of \eqref{1.6}, if and only if
\begin{equation*}
Q=\mathcal{R}_\tau(Q^{p-1})\text{ \ in \ }\R^d.
\end{equation*}

We first establish the following operator-norm continuity of $\mathcal{R}_\tau:L^{p'}_{\rm rad}(\R^d)\to L^p_{\rm rad}(\R^d)$ with respect to $\tau\in[0,1]$.

\begin{Proposition} \label{prop:resolvent}
Under the assumption \eqref{eq:admis}, we have
\begin{equation}\label{eq:reonorm}
\|\mathcal{R}_\tau-\mathcal{R}_{\tau_0}\|_{L^{p'}_{\rm rad}\to L^p_{\rm rad}}\to0\ \ \text{as}\ \ \tau\to \tau_0\in[0,1],
\end{equation}
where $\mathcal{R}_\tau$ is as in \eqref{cont:B-R-tau} for $0\leq \tau\leq1$.
\end{Proposition}

\noindent\textbf{Proof.}
We first consider the endpoint case $\tau_0=1$.
Let $\tau_n\in[0,1)$ satisfy $\tau_n\nearrow1$ as $n\to\infty$. We first claim that $\mathfrak h_{s,\tau_n}\to \mathfrak h_{s,1}$ on $L^2(\R^d)$ as $n\to\infty$ in the sense of \cite[Definition 2.1.1]{Mosco1994}. Actually, we regard $\mathfrak h_{s,\tau}$ as an extended-valued functional on $L^2(\R^d)$, where we define $\mathfrak h_{s,\tau}[u]=+\infty$ for $u\notin D(\mathfrak h_{s,\tau})$. We then have for every $v\in H^s(\R^d)$,
\begin{equation}\label{eq:formmono}
\mathfrak h_{s,\tau_n}[v]=\mathfrak h_{s,1}[v]+(1-\tau_n)\Cc_{d,s}\int_{\R^d}|x|^{-2s}|v|^2 dx\geq\mathfrak h_{s,1}[v].
\end{equation}
Let $u_n\rightharpoonup u$ in $L^2(\R^d)$ as $n\to\infty$. If $\liminf_{n\to\infty}\mathfrak h_{s,\tau_n}[u_n]=+\infty$, then we immediately have $\mathfrak h_{s,1}[u]\leq\liminf_{n\to\infty}\mathfrak h_{s,\tau_n}[u_n].$
If $\liminf_{n\to\infty}\mathfrak h_{s,\tau_n}[u_n]<+\infty$, then we derive from \eqref{eq:formmono} that up to a subsequence if necessary,
\begin{equation*}
 u_n\rightharpoonup u \text{ \ weakly in \ }\mathcal Q^s\text{ \ as \ }n\to\infty.
\end{equation*}
Together with the weakly lower semicontinuity, these imply from \eqref{eq:formmono} that
\begin{equation*}
\mathfrak h_{s,1}[u]\leq\liminf_{n\to\infty}\mathfrak h_{s,1}[u_n]\leq\liminf_{n\to\infty}\mathfrak h_{s,\tau_n}[u_n].
\end{equation*}
On the other hand, for $u\in\mathcal Q^s$, choose $\phi_k\in C_c^\infty(\R^d\setminus\{0\})$ such that
$\phi_k\to u$ strongly in $\mathcal Q^s$ as $k\to\infty$. For every fixed $k$, it follows from \eqref{eq:formmono} that $\mathfrak h_{s,\tau_n}[\phi_k]\to\mathfrak h_{s,1}[\phi_k]$ as $n\to\infty$. By a diagonal choice, then there exist positive integers $k_n\to\infty$ as $n\to\infty$ such that for all sufficiently large $n>0$,
\begin{equation*}
0\leq\mathfrak h_{s,\tau_n}[\phi_{k_n}]-\mathfrak h_{s,1}[\phi_{k_n}]=(1-\tau_n)\Cc_{d,s}\int_{\R^d}|x|^{-2s}|\phi_{k_n}|^2 dx\leq\frac1{k_n}.
\end{equation*}
Since $\phi_k\to u$ strongly in $\mathcal Q^s$ as $k\to\infty$, it gives that
\begin{equation*}
\phi_{k_n}\to u\text{ \ strongly in \  }L^2(\R^d)\ \ \text{as}\ \ n\to\infty,\ \ \mathfrak h_{s,\tau_n}[\phi_{k_n}]\to\mathfrak h_{s,1}[u]\ \ \text{as}\ \ n\to\infty.
\end{equation*}
This proves the claim that $\mathfrak h_{s,\tau_n}\to \mathfrak h_{s,1}$ on $L^2(\R^d)$  as $n\to\infty$ in the sense of \cite[Definition 2.1.1]{Mosco1994}.

By \cite[Theorem~2.4.1]{Mosco1994} and the uniform lower bound
$A_{s,\tau}+1\geq1$, we then deduce that for every $f\in L^2(\R^d)$,
\begin{equation}\label{eq:resolventL2}
\mathcal{R}_{\tau_n}f\to \mathcal{R}_1f\  \ \text{in}\ \ L^2(\R^d) \ \  \text{as}\ \ n\to\infty.
\end{equation}
It follows from \eqref{2.4} that for $0\leq\tau\leq1$,
\begin{equation}\label{eq:unifLp}
\|u\|_p^2\leq C\mathfrak h_{s,1}[u]\leq C\mathfrak h_{s,\tau}[u],\ \  u\in D(\mathfrak h_{s,\tau}).
\end{equation}
Denote $u:=\mathcal{R}_\tau f$ for $f\in L^{p'}_{\rm rad}(\R^d)$. We then obtain from \eqref{eq:varreso} and \eqref{eq:unifLp} that
\begin{equation*}
\mathfrak h_{s,\tau}[u]=\left|\int_{\R^d}f u dx\right|\leq\|f\|_{p'}\|u\|_p\leq C\|f\|_{p'}\mathfrak h_{s,\tau}[u]^{1/2}.
\end{equation*}
This then implies from \eqref{eq:unifLp} that
\begin{equation}\label{eq:resolventb}
\mathfrak h_{s,\tau}[\mathcal{R}_\tau f]^{1/2}+\|\mathcal{R}_\tau f\|_p\leq C\|f\|_{p'},\ \  0\leq\tau\leq1.
\end{equation}
Since $\mathfrak h_{s,1}[\mathcal{R}_\tau f]\leq\mathfrak h_{s,\tau}[\mathcal{R}_\tau f]$, it follows from \eqref{2.7} that
\begin{equation*}
\bigl\{\mathcal{R}_\tau f:~0\leq\tau\leq1,\quad f\in L^{p'}_{\rm rad}(\R^d),\quad \|f\|_{p'}\leq1\bigr\}
\end{equation*}
is relatively compact in $L^p(\R^d)$. Together with \eqref{eq:resolventL2}, this implies that $\mathcal{R}_{\tau_n}g\to \mathcal{R}_1g$ strongly in $L^p(\R^d)$ as $n\to\infty$ for every $g\in L^2_{\rm rad}(\R^d)\cap L^{p'}_{\rm rad}(\R^d)$. It then implies from \eqref{eq:resolventb} that for every $g\in L^{p'}_{\rm rad}(\R^d)$,
\begin{equation}\label{eq:pointLp}
\mathcal{R}_{\tau_n}g\to \mathcal{R}_1g\ \ \ \ \text{in}\  \ L^p(\R^d),
\end{equation}
where we  have used the density of $L^2_{\rm rad}(\R^d)\cap L^{p'}_{\rm rad}(\R^d)$ in $L^{p'}_{\rm rad}(\R^d)$. Moreover, we also have
\begin{equation}\label{eq:selfadjoint}
\int_{\R^d}(\mathcal{R}_\tau f) g dx=\int_{\R^d}f{\mathcal{R}_\tau g} dx,\ \ \forall~ f,g\in L^{p'}_{\rm rad}(\R^d).
\end{equation}

On the contrary, suppose \eqref{eq:reonorm} is false. After passing to a subsequence, then there exist $f_n\in L^{p'}_{\rm rad}(\R^d)$ satisfying $\|f_n\|_{p'}\leq1$ and $\varepsilon>0$ such that for all $n\in\N^+$,
\begin{equation}\label{eq:normcontr}
\|(\mathcal{R}_{\tau_n}-\mathcal{R}_1)f_n\|_p\geq\varepsilon.
\end{equation}
By \eqref{eq:resolventb} and the compact embedding \eqref{2.7}, after passing to a further subsequence, there exists
$0\not\equiv w\in L^p_{\rm rad}(\R^d)$ such that
\begin{equation*}
(\mathcal{R}_{\tau_n}-\mathcal{R}_1)f_n\to w\ \ \text{strongly in}\ \ L^p_{\rm rad}(\R^d)\ \ \text{as}\ \ n\to\infty.
\end{equation*}
We then obtain from \eqref{eq:pointLp} and \eqref{eq:selfadjoint} that for every $g\in L^{p'}_{\rm rad}(\R^d)$,
\begin{equation*}
\begin{aligned}
\int_{\R^d}w g dx
&=\lim_{n\to\infty}\int_{\R^d}(\mathcal{R}_{\tau_n}-\mathcal{R}_1)f_n g dx
 =\lim_{n\to\infty}\int_{\R^d}f_n{(\mathcal{R}_{\tau_n}-\mathcal{R}_1)g} dx=0.
\end{aligned}
\end{equation*}
Hence $w\equiv0$, which however contradicts \eqref{eq:normcontr}. This proves \eqref{eq:reonorm}. One can use the similar argument to deal with the case where $\tau_0\in[0,1)$. This thus completes the proof of Proposition \ref{prop:resolvent}.
\qed

We next construct a locally unique branch of ground states for \eqref{1.6} near the critical coupling $\tau=1$. The construction relies on the resolvent continuity and radial nondegeneracy of $L_{+,1}(Q_1)$.
Define
\begin{equation*}
\mathcal F(\tau,u):=u-\mathcal{R}_\tau\bigl((u_+)^{p-1}\bigr),\ \ (\tau,u)\in[0,1]\times L^p_{\rm rad}(\R^d),
\end{equation*}
where $u_+:=\max\{u,0\}$. Since $p>2$, the map $u\mapsto(u_+)^{p-1}$ is $C^1$ from $L^p_{\rm rad}(\R^d)$ to $L^{p'}_{\rm rad}(\R^d)$. By Proposition~\ref{prop:resolvent}, $\mathcal F$ is continuous, and its derivative in $u$ is also continuous in $(\tau,u)$, where
\begin{equation*}
D_u\mathcal F(\tau,u)h=h-(p-1)\mathcal{R}_\tau\bigl((u_+)^{p-2}h\bigr),\ \ h\in L^p_{\rm rad}(\R^d).
\end{equation*}

\begin{Lemma} \label{lem:endIFT}
Let $Q_1>0$ be a ground state of \eqref{1.6} at $\tau=1$.
Then
there exists $\delta\in(0,1)$ such that \eqref{1.6} admits locally a unique  continuous branch
\begin{equation}\label{eq:endpointbranch}
 \tau\in [1-\delta,1]\longmapsto Q_\tau\in L^p_{\rm rad}(\R^d),
\end{equation}
which bifurcates from $\tau=1$.
Moreover, the map
$$\tau\in (1-\delta,1) \longmapsto
Q_\tau\in H^s_{\rm rad}(\R^d)$$
is of class $C^1$.

Furthermore, for any radial ground state  $Q_*$ of \eqref{1.6} at $\tau_*\in[0,1]$, if $Q_\tau$ is a positive radial weak solution  of \eqref{1.6} at $\tau$, where $|\tau-\tau_*|+\|Q_\tau-Q_*\|_p>0$ is  sufficiently small, then $Q_\tau$ is also a ground state of \eqref{1.6} at $\tau$.
\end{Lemma}

\noindent\textbf{Proof.}
Since $Q_1\in L^p_{\rm rad}(\mathbb R^d)$, one can check that
$$
\|Q_1^{p-2}h\|_{L^{p'}(\mathbb R^d)}
\leq \|Q_1\|_{L^p(\mathbb R^d)}^{p-2}
\|h\|_{L^p(\mathbb R^d)},\ \ \ \forall h\in L^p_{\rm rad}(\mathbb R^d).
$$
Thus, the map
$$
h\longmapsto Q_1^{p-2}h
$$
is bounded from $L^p_{\rm rad}(\mathbb R^d)$ to $L^{p'}_{\rm rad}(\mathbb R^d)$.
We obtain from \eqref{eq:resolventb} and \eqref{2.7} that $$\mathcal{R}_1:~L^{p'}_{\rm rad}(\R^d)\to L^p_{\rm rad}(\R^d)$$
is compact. Thus, $D_u\mathcal F(1,Q_1)$ is Fredholm of index zero. Note from \eqref{eq:varreso} that the identity
$h=(p-1)\mathcal{R}_1(Q_1^{p-2}h)$
is equivalent to $h\in\mathcal Q^s_{\rm rad}$ satisfying $L_{+,1}(Q_1)h=0$ in the form sense. By Proposition~\ref{prop:radialnon}, we then have
\begin{equation*}
\ker D_u\mathcal F(1,Q_1)=\ker\bigl(L_{+,1}(Q_1)|_{\rm rad}\bigr)=\{0\}.
\end{equation*}
Consequently, $D_u\mathcal F(1,Q_1)$ is invertible on $L^p_{\rm rad}(\R^d)$. By the implicit function theorem (cf. \cite[Proposition~2.1]{Glockner2006}), we thus obtain that the continuous branch \eqref{eq:endpointbranch} exists and is locally unique in an $L^p_{\rm rad}(\R^d)$-neighborhood of $Q_1$. We then have $Q_\tau\not\equiv0$ by decreasing $\delta>0$ if necessary. Since
$Q_\tau=\mathcal{R}_\tau\bigl(((Q_\tau)_+)^{p-1}\bigr)$,
we deduce from \cite[Theorem~1.1]{BogdanEtAl2019} that $Q_\tau>0$. Moreover, $D_u\mathcal F(\tau,Q_\tau)$ remains invertible on $L^p_{\rm rad}(\R^d)$ for every $\tau\in[1-\delta,1]$.

For $0<\tau<1$, we have $Q_\tau\in H^s_{\rm rad}(\R^d)$. Consider the $C^1$ map
\begin{equation*}
\mathcal G:(0,1)\times H^s_{\rm rad}(\R^d)\longrightarrow H^{-s}_{\rm rad}(\R^d),
\end{equation*}
which is defined by
\begin{equation*}
\langle\mathcal G(\tau,u),v\rangle:=\mathfrak h_{s,\tau}[u,v]-\int_{\R^d}(u_+)^{p-1}v dx,\ \ u,v\in H^s_{\rm rad}(\R^d).
\end{equation*}
One can check that $\ker D_u\mathcal G(\tau,Q_\tau)=\{0\}$ and $D_u\mathcal G(\tau,Q_\tau)$ is Fredholm of index zero. Therefore, $D_u\mathcal G(\tau,Q_\tau)$ is an isomorphism. The implicit function theorem then implies that \eqref{1.6} admits locally a unique $C^1-$branch  in $L^p_{\rm rad}(\R^d)\cap H^s_{\rm rad}(\R^d)$. The local uniqueness in $L^p_{\rm rad}(\R^d)$ then identifies it with \eqref{eq:endpointbranch}.

We next prove that Morse index one is locally preserved in the positive radial solution set. Set for $0\leq\tau\leq1$ and $0\leq Q\in L^p_{\rm rad}(\R^d)$,
\begin{equation*}
m_Q:=\sqrt{p-1}Q^{(p-2)/2},\ \ r:=\frac{2p}{p-2},
\end{equation*}
and define
\begin{equation}\label{eq:BSope}
\mathbb K_{\tau,Q}f:=m_Q\mathcal{R}_\tau(m_Qf),\ \ f\in L^2_{\rm rad}(\R^d).
\end{equation}
Since $m_Q\in L^r_{\rm rad}(\R^d)$, we obtain from H\"older's inequality that
\begin{equation*}
\|m_Qf\|_{p'}\leq\|m_Q\|_r\|f\|_2,\ \ \|m_Qu\|_2\leq\|m_Q\|_r\|u\|_p.
\end{equation*}
This then implies that
\begin{equation*}
\begin{aligned}
\|\mathbb K_{\tau,Q}f\|_{L^2_{\rm rad}}
&=\|m_Q\mathcal R_\tau(m_Qf)\|_{L^2_{\rm rad}}
\leq \|m_Q\|_{L^r_{\rm rad}}
\|\mathcal R_\tau(m_Qf)\|_{L^p_{\rm rad}}\\
&\leq \|m_Q\|_{L^r_{\rm rad}}
\|\mathcal R_\tau\|_{\mathcal L(L^{p'}_{\rm rad},L^p_{\rm rad})}
\|m_Qf\|_{L^{p'}_{\rm rad}}\\
&\leq \|m_Q\|_{L^r_{\rm rad}}^2
\|\mathcal R_\tau\|_{\mathcal L(L^{p'}_{\rm rad},L^p_{\rm rad})}
\|f\|_{L^2_{\rm rad}}.
\end{aligned}
\end{equation*}
Consequently, we have
\begin{equation*}
\|\mathbb K_{\tau,Q}\|_{\mathcal L(L^2_{\rm rad},L^2_{\rm rad})}
\leq
\|m_Q\|_{L^r_{\rm rad}}^2
\|\mathcal R_\tau\|_{\mathcal L(L^{p'}_{\rm rad},L^p_{\rm rad})},
\end{equation*}
where $\mathcal L(A,B)$ denotes the space of bounded linear operators
from a Banach space $A$ to a Banach space $B$, equipped with the operator norm. Since $\mathcal{R}_\tau:L^{p'}_{\rm rad}(\R^d)\to L^p_{\rm rad}(\R^d)$ is compact, $\mathbb K_{\tau,Q}$ is compact on $L^2_{\rm rad}(\R^d)$.

It follows from \eqref{eq:varreso} that for any $f,g\in L^2_{\rm rad}(\R^d)$,
\begin{equation*}
\begin{aligned}
(\mathbb K_{\tau,Q}f,g)_{L^2}=&\mathfrak h_{s,\tau}[\mathcal{R}_\tau(m_Qf),\mathcal{R}_\tau(m_Qg)]
= (f,\mathbb K_{\tau,Q}g)_{L^2}.
\end{aligned}
\end{equation*}
Moreover, we have
\begin{equation*}
(\mathbb K_{\tau,Q}f,f)_{L^2}=\mathfrak h_{s,\tau}[\mathcal{R}_\tau(m_Qf)]\geq0.
\end{equation*}
Thus, $\mathbb K_{\tau,Q}$ is nonnegative and self-adjoint.
The Birman--Schwinger principle of closed forms (cf. \cite[Proposition~2.1]{LaptevSolomyak2012}) yields that
\begin{equation}\label{eq:BSindex}
\operatorname{ind}_{\rm rad}(L_{+,\tau}(Q))=\#\{j:~\lambda_j(\mathbb K_{\tau,Q})>1\},
\end{equation}
where the eigenvalues $\lambda_j(\mathbb K_{\tau,Q})$ are counted with multiplicity.
Here and below, $\operatorname{ind}_{\rm rad}(L_{+,\tau}(Q))$
denotes the Morse index of $L_{+,\tau}(Q)$ restricted to radial functions. Moreover, we get from \eqref{eq:varreso} that
\begin{equation}\label{eq:BSkernel}
\dim\ker\bigl(L_{+,\tau}(Q)|_{\rm rad}\bigr)=\dim\ker(I-\mathbb K_{\tau,Q}).
\end{equation}

Suppose $\tau_n\in[0,1]$ satisfy $\tau_n\to\tau$ as $n\to\infty$, and let $Q_{\tau_n}\geq0$ satisfy
\begin{equation*}
Q_{\tau_n}\to Q_{\tau}\text{ \ strongly in \  }L^p_{\rm rad}(\R^d)\text{ \ as \  }n\to\infty,
\end{equation*}
where $Q_\tau>0$ is a positive radial weak solution  of \eqref{1.6} at $\tau$.
This then implies that for $r=\frac{2p}{p-2}$,
\begin{equation*}
m_{Q_{\tau_n}}\to m_{Q_{\tau}}\text{ \ strongly in \  }L^r(\R^d)\text{ \ as \  }n\to\infty.
\end{equation*}
By Proposition~\ref{prop:resolvent}, we obtain from \eqref{eq:resolventb} and \eqref{eq:BSope} that
\begin{equation}\label{eq:BSnorm}
\begin{aligned}
&\|\mathbb K_{\tau_n,Q_{\tau_n}}-\mathbb K_{\tau,Q_{\tau}}\|_{\mathcal L(L^2_{\rm rad},L^2_{\rm rad})}\\
\leq&C\bigl(\|m_{Q_{\tau_n}}\|_r+\|m_{Q_{\tau}}\|_r\bigr)\|m_{Q_{\tau_n}}-m_{Q_{\tau}}\|_r+\|m_{Q_{\tau}}\|_r^2\|\mathcal{R}_{\tau_n}-\mathcal{R}_\tau\|_{\mathcal L(L^{p'}_{\rm rad},L^p_{\rm rad})}\\
\to&0\text{ \ as \  }n\to\infty.
\end{aligned}
\end{equation}
Note from Proposition~\ref{prop:radialnon} that $Q_*>0$ is radially nondegenerate. Since $Q_*$ is a ground state of \eqref{1.6} with $\tau=\tau_*$, we obtain from \eqref{eq:BSindex}, \eqref{eq:BSkernel} and the compactness of $\mathbb K_{\tau_*,Q_*}$ that $\mathbb K_{\tau_*,Q_*}$ has exactly one eigenvalue in $(1,\infty)$, counted with multiplicity, and $1\notin\sigma(\mathbb K_{\tau_*,Q_*})$. Thus, $\mathbb K_{\tau_*,Q_*}$ has a spectral gap around $1$. We then derive from \eqref{eq:BSnorm} and \cite[Chapter~IV, Theorem~3.16]{Kato1995} that $\mathbb K_{\tau,Q_{\tau}}$ also has exactly one eigenvalue in $(1,\infty)$, counted with multiplicity, where $Q_\tau$ is a positive radial weak
solution of \eqref{1.6} at $\tau$ and $|\tau-\tau_*|+\|Q_{\tau}-Q_*\|_p>0$ is sufficiently small. Together with \eqref{eq:BSindex}, this shows that these solutions
have radial Morse index one. Applying Proposition~\ref{prop:1channel},  the angular comparison used in \eqref{eq:fullcom} then gives that their nonradial channels are nonnegative, and hence their full
Morse indexes are one. This completes the proof of Lemma~\ref{lem:endIFT}.
\qed
\vspace{0.5em}

We next establish the strong compactness along ground state branches away from the critical coupling $\tau=1$.

\begin{Proposition} \label{prop:strongcon}
Let $\tau\mapsto Q_\tau$ belong to $C^1(I,H^s_{\rm rad}(\R^d))$, where $I\subset(0,1)$ is an open interval and $Q_\tau$ is a ground state for \eqref{1.6}.
Suppose that $\{\tau_n\}\subset I$ satisfies $0<\tau_n\leq\tau_0<1$ and $\tau_n\to\tau_*\in[0,\tau_0]$ as $n\to\infty$. Then there exists $Q_*\in H^s_{\rm rad}(\R^d)$ such that up to a subsequence,
\begin{equation}\label{eq:strongHs}
Q_{\tau_n}\to Q_*\text{ \ strongly in \  }H^s_{\rm rad}(\R^d)\text{ \ as \  }n\to\infty,
\end{equation}
where $Q_*>0$ is a ground state of \eqref{1.6} at $\tau=\tau_*$ and is radially nondegenerate.
\end{Proposition}

\noindent\textbf{Proof.}
For $0<\tau<1$ in the branch interval, set
\begin{equation*}
f(\tau):=\int_{\R^d}Q_\tau^p dx.
\end{equation*}
We then obtain from \eqref{1.6} that
\begin{equation}\label{eq:nehari}
\mathfrak h_{s,\tau}[Q_\tau]=\mathfrak q_{s,\tau}[Q_\tau]+\|Q_\tau\|_2^2=f(\tau).
\end{equation}
Since $\tau\mapsto Q_\tau$ belongs to
$C^1(I, H^s_{\rm rad}(\R^d))$, differentiating \eqref{eq:nehari} and multiplying \eqref{1.6} by $\partial_\tau Q_\tau$, it yields that
\begin{equation*}
\begin{aligned}
f'(\tau)=&-\Cc_{d,s}\int_{\R^d}|x|^{-2s}Q_\tau^2 dx+2\mathfrak h_{s,\tau}[Q_\tau,\partial_\tau Q_\tau]\\
=&-\Cc_{d,s}\int_{\R^d}|x|^{-2s}Q_\tau^2 dx+2\int_{\R^d}Q_\tau^{p-1}\partial_\tau Q_\tau dx\\
=&-\Cc_{d,s}\int_{\R^d}|x|^{-2s}Q_\tau^2 dx+\frac2p f'(\tau),
\end{aligned}
\end{equation*}
which implies that
\begin{equation}\label{cont:V-prime}
f'(\tau)=-\frac{p\Cc_{d,s}}{p-2}\int_{\R^d}|x|^{-2s}Q_\tau^2 dx<0.
\end{equation}

By Hardy's inequality, we obtain from \eqref{eq:nehari} that
\begin{equation*}
(1-\tau)\Cc_{d,s}\int_{\R^d}|x|^{-2s}Q_\tau^2 dx\leq\mathfrak q_{s,\tau}[Q_\tau]\leq f(\tau).
\end{equation*}
Together with \eqref{cont:V-prime}, this implies that
$0<-\frac{f'(\tau)}{f(\tau)}\leq\frac{p}{(p-2)(1-\tau)}$.
For $0<\tau\leq\widehat\tau<1$ in the branch interval, integrating over $[\tau,\widehat\tau]$, we deduce that
$\log\frac{f(\tau)}{f(\widehat\tau)}\leq\frac{p}{p-2}\log\frac{1-\tau}{1-\widehat\tau}$.
This implies that
\begin{equation}\label{eq:fupperb}
f(\tau)\leq f(\widehat\tau)\left(\frac{1-\tau}{1-\widehat\tau}\right)^{p/(p-2)}.
\end{equation}
On the other hand, it follows from \eqref{eq:unifLp} and \eqref{eq:nehari} that
\begin{equation*}
f(\tau)^{2/p}=\|Q_\tau\|_p^2\leq C\mathfrak h_{s,\tau}[Q_\tau]=Cf(\tau).
\end{equation*}
Since $f(\tau)>0$ and $p>2$, we have
\begin{equation}\label{eq:flowerb}
f(\tau)\geq C^{-p/(p-2)}=:c_{d,s,p}>0.
\end{equation}

Write $Q_n:=Q_{\tau_n}$ and fix $\widehat\tau\in I$. We obtain from \eqref{cont:V-prime} and \eqref{eq:fupperb} that
\begin{equation*}
f(\tau_n)\leq f(\widehat\tau)(1-\widehat\tau)^{-p/(p-2)},\ \ n\geq1.
\end{equation*}
By Hardy's inequality, it also implies that
\begin{equation}\label{eq:unicoer}
\mathfrak h_{s,\tau_n}[u]\geq(1-\tau_0)\|(-\Delta)^{s/2}u\|_2^2+\|u\|_2^2\geq c_{\tau_0}\|u\|_{H^s}^2,\ \ u\in H^s(\R^d),
\end{equation}
where $c_{\tau_0}>0$ is independent of $n\geq1$. Together with \eqref{eq:nehari}, this implies that $\{Q_n\}$ is bounded uniformly in $H^s(\R^d)$. We thus get from \eqref{2.7} that up to a subsequence, there exists $Q_*\in H^s_{\rm rad}(\R^d)$ such that
\begin{equation}\label{eq:precompact}
\begin{aligned}
&Q_n\rightharpoonup Q_*\text{ \ weakly in \  }H^s_{\rm rad}(\R^d)\ \ \text{as}\ \ n\to\infty, \\
&Q_n\to Q_*\text{ \ strongly in \  }L^p_{\rm rad}(\R^d)\ \ \text{as}\ \ n\to\infty.
\end{aligned}
\end{equation}
Since $Q_n>0$, we have $Q_*\geq0$ $a.e.$ in $\R^d$. Passing to the weak limit in \eqref{1.6}, it yields from \eqref{eq:precompact} that
\begin{equation*}
(A_{s,\tau_*}+1)Q_*=Q_*^{p-1}\text{ \ in \  }\R^d.
\end{equation*}
It then yields from \eqref{eq:flowerb} and \eqref{eq:precompact} that
\begin{equation*}
\|Q_*\|_p^p=\lim_{n\to\infty}f(\tau_n)\geq c_{d,s,p}>0.
\end{equation*}
Together with \cite[Theorem~1.1 and Section~2.1]{BogdanEtAl2019}, we conclude from \eqref{eq:varreso} that $Q_*>0$.
Let $w_n:=Q_n-Q_*$. Using the equations of $Q_n$ and $Q_*$, we get that
\begin{equation}\label{eq:denergy}
\mathfrak h_{s,\tau_n}[w_n]=\int_{\R^d}(Q_n^{p-1}-Q_*^{p-1})w_n dx+(\tau_n-\tau_*)\Cc_{d,s}\int_{\R^d}|x|^{-2s}Q_*w_n dx.
\end{equation}
Applying H\"older's inequality and Hardy's inequality, we deduce from \eqref{eq:unicoer} and \eqref{eq:denergy} that
\begin{equation*}
c_{\tau_0}\|w_n\|_{H^s}^2\leq\mathfrak h_{s,\tau_n}[w_n]\leq C\left(\|Q_n^{p-1}-Q_*^{p-1}\|_{p'}+|\tau_n-\tau_*|\|Q_*\|_{H^s}\right)\|w_n\|_{H^s}.
\end{equation*}
Together with \eqref{eq:precompact}, this proves \eqref{eq:strongHs}.

If $\mathbb K_{\tau_*,Q_*}$ had at least two eigenvalues in $(1,\infty)$, counted with multiplicity, it would follow from \eqref{eq:BSnorm} that $\mathbb K_{\tau_n,Q_n}$ also has at least two eigenvalues in $(1,\infty)$ for sufficiently large $n$. This however contradicts \eqref{eq:BSindex}, since any $Q_n$ is a ground state of \eqref{1.6} with $\tau=\tau_n$. We thus deduce from \eqref{eq:BSindex} that
\begin{equation*}
\operatorname{ind}_{\rm rad}(L_{+,\tau_*}(Q_*))\leq1.
\end{equation*}
On the other hand, the equation of $Q_*$ implies that
\begin{equation*}
\mathfrak h_{s,\tau_*}[Q_*]-(p-1)\int_{\R^d}Q_*^p dx=-(p-2)\int_{\R^d}Q_*^p dx<0.
\end{equation*}
By the angular comparison used in \eqref{eq:fullcom}, we thus conclude from Proposition~\ref{prop:1channel} that
\begin{equation*}
\operatorname{ind}(L_{+,\tau_*}(Q_*)) =\operatorname{ind}_{\rm rad}(L_{+,\tau_*}(Q_*))=1.
\end{equation*}
Proposition~\ref{prop:radialnon} then implies that $Q_*$ is radially nondegenerate. This completes the proof of Proposition~\ref{prop:strongcon}.\qed
\vspace{0.5em}

\subsection{Proof of Theorem~\ref{thm:main} }\label{sub3.1}
Combining the local continuation with the strong compactness of Proposition~\ref{prop:strongcon}, in this subsection we now complete the proof of Theorem~\ref{thm:main}.

\noindent\textbf{Proof of Theorem~\ref{thm:main}.}
Let $Q_1$ be a  ground state of \eqref{1.6} at $\tau=1$. We first show that the local branch obtained in Lemma~\ref{lem:endIFT} can be extended uniquely to a continuous branch of ground states
\begin{equation*}
\tau\in[0,1]\longmapsto Q_\tau\in L^p_{\rm rad}(\R^d),
\end{equation*}
which is continuous in $H^s_{\rm rad}(\R^d)$ on $[0,1)$ and is $C^1$ in $\tau\in (0,1)$.
By Proposition~\ref{prop:radialnon} and Lemma~\ref{lem:endIFT}, the implicit function theorem applies at every point of the branch with $0<\tau<1$. Hence, the branch extends locally and uniquely through any point of this type, while it preserves the positivity and Morse index one.
Let $(a,1]$, where $0\leq a<1$, be the maximal
interval of this extension within $(0,1]$.
The map $\tau\mapsto Q_\tau$ is of class
$C^1$ from $(a,1)$ into $H^s_{\rm rad}(\R^d)$.

We now claim that $a=0$. On the contrary, suppose that $a>0$. By Proposition~\ref{prop:strongcon}, there exist a sequence $\tau_n\searrow a$ as $n\to\infty$ and a ground state $Q_a$ of \eqref{1.6} at $\tau=a$ such that
\begin{equation}\label{eq:6.35}
Q_{\tau_n}\to Q_a\text{ \ strongly in \  }H^s_{\rm rad}(\R^d)\text{ \ as \  }n\to\infty.
\end{equation}
By Proposition~\ref{prop:radialnon}, $Q_a$ is radially nondegenerate. Arguing as in Lemma~\ref{lem:endIFT}, we then deduce that there exists a locally unique branch $\widetilde Q_\tau$ through $(a,Q_a)$, defined for $\tau$ on both sides of $a$. Together with \eqref{eq:6.35}, this implies that $Q_{\tau_n}\equiv\widetilde Q_{\tau_n}$ in $\R^d\setminus\{0\}$ for all sufficiently large $n$. By the local uniqueness of the positive radial solution
branch through $(a,Q_a)$, we thus obtain that there exists a sufficiently small $\eta>0$ such that
\begin{equation*}
Q_\tau\equiv\widetilde Q_\tau\text{ \ in \  }\R^d\setminus\{0\},\ \ a<\tau<a+\eta.
\end{equation*}
Moreover, we derive from Lemma~\ref{lem:endIFT} that $\widetilde Q_\tau$ is a positive ground state of \eqref{1.6} for $\tau<a$ sufficiently close to $a$. This extends the original branch below $a$, which contradicts the maximality of $(a,1]$. We therefore conclude that $a=0$.

Finally, choose $\tau_n\in(0,1)$ such that $\tau_n\searrow0$ as $n\to\infty$. By Proposition~\ref{prop:strongcon}, after passing to a subsequence, there exists a ground state $Q_0$ of \eqref{1.6} at $\tau=0$ such that $Q_{\tau_n}\to Q_0$ strongly in $H^s_{\rm rad}(\R^d)$ as $n\to\infty$. By the radial nondegeneracy of $Q_0$ and the implicit function theorem (cf. \cite[Proposition~2.1]{Glockner2006}), we then get that
\begin{equation*}
Q_\tau\to Q_0\text{ \ strongly in \  }H^s_{\rm rad}(\R^d)\text{ \ as \  }\tau\searrow0.
\end{equation*}
This proves the claimed unique continuous extension to $\tau=0$.
By Theorem~\ref{thm1.1}, there exists a radial ground state of \eqref{1.6} at $\tau=1$. Let $Q_1^{(1)}$ and $Q_1^{(2)}$ be two ground states of \eqref{1.6} at $\tau=1$. Applying the continuation of each ground state, we derive that $Q_1^{(1)}$ and $Q_1^{(2)}$ generate two continuous branches
$$(Q_\tau^{(i)})_{\tau\in[0,1]}
\in C\big([0,1], L^p_{\rm rad}(\R^d)\big),
\qquad i=1,2.$$
The endpoints $Q_0^{(1)},Q_0^{(2)}\in H^s_{\rm rad}(\R^d)$
are ground states of \eqref{1.6} with $\tau=0$. Since $Q_0^{(1)},Q_0^{(2)}\in H^s_{\rm rad}$, we obtain from \cite[Theorem~3.4]{FLS2016} that
\begin{equation}\label{eq:tau0}
Q_0^{(1)}\equiv Q_0^{(2)}\text{ \ in \  }\R^d.
\end{equation}

Define
\begin{equation*}
E:=\left\{\tau\in[0,1]:~\|Q_\tau^{(1)}-Q_\tau^{(2)}\|_p=0\right\}.
\end{equation*}
It then follows from \eqref{eq:tau0} that $0\in E$, and $E$ is nonempty. By the continuity of both branches in $L^p_{\rm rad}(\R^d)$, $E$ is relatively closed in $[0,1]$. Set $Q_*:=Q^{(1)}_{\tau_*}\equiv Q^{(2)}_{\tau_*}$ for every $\tau_*\in E$. By Proposition~\ref{prop:radialnon}, $Q_*$ is radially nondegenerate. Applying the implicit function theorem, we thus obtain that there exists a unique local radial solution branch through $(\tau_*,Q_*)$. This further implies that for sufficiently small $\eta>0$,
\begin{equation*}
Q_\tau^{(1)}\equiv Q_\tau^{(2)}\text{ \ in \  }\R^d\setminus\{0\},\ \ \tau\in[0,1]\cap(\tau_*-\eta,\tau_*+\eta).
\end{equation*}
Thus, $E$ is relatively open in $[0,1]$. Since $[0,1]$ is connected, we conclude that $E=[0,1]$. In particular, we get that
\begin{equation*}
Q_\tau^{(1)}\equiv Q_\tau^{(2)}\text{ \ in \  }\R^d\setminus\{0\}.
\end{equation*}
This proves the uniqueness of ground states for \eqref{1.6} at $\tau=1$. For any $\tau_*\in(0,1)$, the same continuation argument extends every ground state from $\tau=\tau_*$ to $\tau=0$. By \cite[Theorem~3.4]{FLS2016} and the local uniqueness along the branches, the resulting branch coincides with $Q_\tau^{(1)}$ on $[0,\tau_*]$. This proves the uniqueness of ground states of \eqref{1.6} for every $0<\tau\leq1$.

For every $0<\tau\leq1$, let $Q_\tau$ denote this unique ground state of \eqref{1.6}. We then obtain from Theorem~\ref{thm:fullkernel} that $\ker L_{+,\tau}(Q_\tau)=\{0\}$.
This completes the proof of Theorem~\ref{thm:main}.
\qed

\subsection{Applications of Theorem~\ref{thm:main}}
\underline{}\label{sec:completion}
Applying Theorem~\ref{thm:main}, in this subsection we mainly establish Theorem~\ref{thm:sharp-HGN} and Corollary \ref{cor:one-particle-HLT}. We start with the proof of Theorem~\ref{thm:sharp-HGN}.
\vspace{0.5em}

\noindent\textbf{Proof of Theorem~\ref{thm:sharp-HGN}.}
Since the existence of Theorem~\ref{thm:sharp-HGN}~(1) is proved in Appendix~\ref{sec:forms},
here it suffices to prove Theorem~\ref{thm:sharp-HGN} (2).

Note that every function in \eqref{1.21} is an optimizer of \eqref{1.17}. Conversely, let $u\in\mathcal Q^s\setminus\{0\}$ be an optimizer of \eqref{1.17}. By the ground-state representation of  \cite[Theorem 1.2]{FrankSeiringer2008}, we obtain that
\begin{equation*}
|u|\in\mathcal Q^s,\ \ \mathfrak q_{s,1}[|u|]\leq\mathfrak q_{s,1}[u].
\end{equation*}
Applying the optimality of $u$, this then implies that
\begin{equation}\label{eq:hgnoptimizer}
\mathfrak q_{s,1}[|u|]=\mathfrak q_{s,1}[u].
\end{equation}
Using the scaling in \eqref{2.13} and \eqref{2.14}, $|u|$ gives a nonnegative minimizer of \eqref{1.7}. By Proposition~\ref{prop2.2} and Theorem~\ref{thm:main}, we then have
\begin{equation}\label{eq:oporbit}
|u(x)|\equiv a\lambda^{d/2}Q_1(\lambda x)\ \ \text{in}\ \ \R^d\setminus\{0\}\text{ \ for some \ }a,\lambda>0.
\end{equation}
Moreover, using the ground state representations of $u$ and $|u|$, we derive from \eqref{eq:hgnoptimizer} that
\begin{equation*}
\iint_{\R^d\times\R^d}\frac{|u(x)||u(y)|-u(x)u(y)}{|x-y|^{d+2s}} dx dy=0.
\end{equation*}
Since the integrand is nonnegative and $|u|>0$ in $\R^d\setminus\{0\}$, we conclude that $u$ does not change the sign in $\R^d\setminus\{0\}$. Together with \eqref{eq:oporbit}, this implies that for some $a,\lambda>0$,
\begin{equation*}
u(x)\equiv\pm a\lambda^{d/2}Q_1(\lambda x)\text{ \ in \  }\R^d\setminus\{0\}.
\end{equation*}
This proves Theorem~\ref{thm:sharp-HGN} (2), which thus completes the proof of Theorem~\ref{thm:sharp-HGN}.
\qed

\vspace{0.5em}
Applying Theorem~\ref{thm:sharp-HGN}, we now prove Corollary~\ref{cor:one-particle-HLT} as follows.
\vspace{0.5em}

\noindent\textbf{Proof of Corollary~\ref{cor:one-particle-HLT}.}
For every $\gamma>0$, set $\bar{r}:=\gamma+\frac{d}{2s}$.
Let $V\in L^{\bar{r}}(\R^d)$ be a real-valued function.
Applying H\"older's inequality, we then obtain from \eqref{1.17} that
\begin{equation*}
\int_{\R^d}Vu^2 dx\leq \|V_+\|_{\bar{r}}\|u\|_p^2\leq C_{\rm HGN}\|V_+\|_{\bar{r}}\mathfrak q_{s,1}[u]^\theta\|u\|_2^{2(1-\theta)},\ \  u\in\mathcal Q^s,
\end{equation*}
where $p=\frac{2\bar{r}}{\bar{r}-1}=\frac{2(d+2s\gamma)}{d+2s\gamma-2s}\in(2,2_s^*)$ and $0<\theta<1$ is as in \eqref{01}.
For $\|u\|_2=1$, minimizing the function $t-C_{\rm HGN}\|V_+\|_{\bar{r}} t^\theta$ over $t\geq0$, we deduce from  \eqref{eq:sharp-HGN-constant-nu} and \eqref{2.10} that
\begin{equation}\label{eq:HLT1}
A_{s,1}-V\geq-\nu_{s,1}^{-1/(1-\theta)}\|V_+\|_{\bar{r}}^{1/(1-\theta)}\ \ \text{on}\ \ L^2(\R^d).
\end{equation}
Since $\gamma={\bar{r}}(1-\theta)$, this proves \eqref{eq:KHLT}.

Define for $\lambda>0$,
\begin{equation*}
V_\lambda(x):=\lambda^{2s}Q_1(\lambda x)^{p-2},\ \ Q^{(\lambda)}(x):=\lambda^{d/2}Q_1(\lambda x)\ \  \ \text{in}\ \ \R^d\setminus\{0\}.
\end{equation*}
We then derive from \eqref{2.10}, \eqref{2.13} and \eqref{2.15} that
\begin{equation*}
\frac{\mathfrak q_{s,1}[Q^{(\lambda)}]-\int_{\R^d}V_\lambda |Q^{(\lambda)}|^2 dx}{\|Q^{(\lambda)}\|_2^2}=-\lambda^{2s},\ \ \int_{\R^d}V_\lambda^{\bar{r}} dx=\lambda^{2s\gamma}\nu_{s,1}^{\bar{r}}.
\end{equation*}
It follows from  \eqref{eq:HLT1} that $\lambda_1(A_{s,1}-V_\lambda)=-\lambda^{2s}$. Thus, the identity of \eqref{eq:KHLT} holds for every $V_\lambda$, where the constant $\nu_{s,1}^{-(\gamma+\frac{d}{2s})}>0$ is optimal.

If $V\in L^{\bar{r}}(\R^d)$ satisfying $\|V_+\|_{\bar{r}}>0$ is an optimizer of \eqref{eq:KHLT}, then we have $\lambda_1(A_{s,1}-V)<0$. Arguing as in \cite[Theorem~11.6]{LiebLoss2001}, we obtain that there exists $u\in\mathcal Q^s$ satisfying $\|u\|_2=1$ such that
\begin{equation}\label{eq:HLT2}
\begin{aligned}
\lambda_1(A_{s,1}-V)
&=\mathfrak q_{s,1}[u]-\int_{\R^d}Vu^2 dx
 \geq \mathfrak q_{s,1}[u]-\|V_+\|_{\bar{r}}\|u\|_p^2\\
&\geq \mathfrak q_{s,1}[u]-C_{\rm HGN}\|V_+\|_{\bar{r}}\mathfrak q_{s,1}[u]^\theta\\
&\geq-\nu_{s,1}^{-1/(1-\theta)}\|V_+\|_{\bar{r}}^{1/(1-\theta)}=\lambda_1(A_{s,1}-V),
\end{aligned}
\end{equation}
where the last inequality follows by the similar proof of \eqref{eq:HLT1}. Thus, all identities of \eqref{eq:HLT2} hold, and hence $u$ is an optimizer of \eqref{1.17}. We then deduce from Theorem~\ref{thm:sharp-HGN} that
\begin{equation*}
|u(x)|\equiv aQ^{(\lambda)}(x)\ \ \text{in}\ \ \R^d\setminus\{0\},~\ \  \forall~a,\lambda>0.
\end{equation*}
We also derive from \eqref{eq:HLT2} that the minimum of the function $t-C_{\rm HGN}\|V_+\|_{\bar{r}} t^\theta$
over $t\geq0$ is attained at $t_0=\mathfrak q_{s,1}[u]$.
Since $\|u\|_2=1$, note from \eqref{1.17} that $\mathfrak q_{s,1}[u]>0$. Thus, the minimizing point $t_0=\mathfrak q_{s,1}[u]$ lies in $(0,\infty)$.
Since the derivative of the function $t-C_{\rm HGN}\|V_+\|_{\bar r}t^\theta$
is zero at $t_0=\mathfrak q_{s,1}[u]$, we obtain that
$$1-\theta C_{\rm HGN}\|V_+\|_{\bar r}\mathfrak q_{s,1}[u]^{\theta-1}=0.$$
By the identity of \eqref{1.17} and $\|u\|_2=1$, this implies that
\begin{equation*}
\mathfrak q_{s,1}[u]=\theta C_{\rm HGN}\|V_+\|_{\bar{r}}\mathfrak q_{s,1}[u]^\theta=\theta\|V_+\|_{\bar{r}}\|u\|_p^2.
\end{equation*}
We therefore obtain from  \eqref{2.10}, \eqref{2.13} and \eqref{2.15} that
\begin{equation*}
\|V_+\|_{\bar{r}}=\frac{\mathfrak q_{s,1}[u]}{\theta\|u\|_p^2}=\lambda^{2s-d/\bar{r}}\|Q_1\|_p^{p-2}.
\end{equation*}
Since $|u|>0$ holds in $\R^d\setminus\{0\}$, if all identities of
\begin{equation*}
\int_{\R^d}Vu^2 dx\leq\int_{\R^d}V_+u^2 dx\leq\|V_+\|_{\bar{r}}\|u\|_p^2
\end{equation*}
hold, then $V_-\equiv 0~a.e.$ in $\R^d\setminus\{0\}$, and
\begin{equation*}
V(x)=\frac{\|V_+\|_{\bar{r}}}{\|u\|_p^{p-2}}|u(x)|^{p-2}=\lambda^{2s}Q_1(\lambda x)^{p-2}\ \ \text{in}\ \ \R^d\setminus\{0\}.
\end{equation*}
This hence proves \eqref{eq:one-particle-HLT-optimizers}. We thus complete the proof of Corollary~\ref{cor:one-particle-HLT}.
\qed
\vspace{0.5em}
\noindent{\bf Acknowledgements:} The authors thank Professor Phan Th\`anh  Nam very much for suggesting this interesting problem and fruitful discussions of this manuscript.
\appendix
\renewcommand{\theequation}{\thesection.\arabic{equation}}
\section{Ground States and Sharp Interpolation Inequality}\label{sec:forms}
This appendix is devoted to  the proofs of Theorems \ref{thm1.1} and \ref{thm:sharp-HGN} (1) on the existence and analytical properties of ground states for \eqref{1.6}.
Let $\mathfrak h_{s,\tau}[\cdot]$ be defined by \eqref{2.2}. We  first prove the following estimates.

\begin{Lemma}\label{lem2.1}
For every $0<t<s<1$, we have
\begin{equation}
\|{u}\|_{H^t}^2\leq C_{d,s,t}\mathfrak h_{s,1}[u], \quad \forall ~ u\in\mathcal Q^s.
\label{2.3}
\end{equation}
For every $2<q<2_s^*$, we also have
\begin{equation}\label{2.4}
 \|u\|_q^2\leq C_{d,s,q}\mathfrak h_{s,\tau}[u]\text{ \ for \  }u\in D(\mathfrak h_{s,\tau}),\ 0\leq\tau\leq1.
\end{equation}
In particular,
\begin{equation}
\mathcal Q^s_{\rm rad}\hookrightarrow L^q(\R^d)\text{ \ is compact for \  }2<q<2_s^*.
\label{2.7}
\end{equation}
\end{Lemma}

\noindent\textbf{Proof.}
By \cite[Theorem~1.2]{Frank2009} and Young's inequality, we have
\begin{equation*}
\begin{aligned}
\|u\|_{H^t}^2\leq&C_{d,s,t}\left(\mathfrak q_{s,1}[u]^{t/s}\|u\|_2^{2(1-t/s)}+\|u\|_2^2\right)\\
\leq&C_{d,s,t}\mathfrak h_{s,1}[u]\text{ \ for \  }u\in C_c^\infty(\R^d\setminus\{0\}),\ 0<t<s.
\end{aligned}
\end{equation*}
By the density of $C_c^\infty(\R^d\setminus\{0\})$ in $\mathcal Q^s$ with respect to the norm $\mathfrak h_{s,1}^{1/2}$, we obtain that \eqref{2.3} holds.

For $2<q<2_s^*$, choose $0<t<s$ such that $q<2_t^*$. By Sobolev's embedding theorem and \eqref{2.3}, we then have
\begin{equation*}
\|u\|_q^2\leq C_{d,s,q}\|u\|_{H^t}^2\leq C_{d,s,q}\mathfrak h_{s,1}[u]\leq C_{d,s,q}\mathfrak h_{s,\tau}[u]\text{ \ for \  }u\in D(\mathfrak h_{s,\tau}),\ 0\leq\tau\leq1,
\end{equation*}
where the last inequality follows from \eqref{2.1} and \eqref{2.2}. Finally, it yields from \eqref{2.3} and \cite[Remark~II.5(i)]{Lions1982}  that \eqref{2.7} holds. This completes the proof of Lemma~\ref{lem2.1}.
\qed
 \vspace{0.5em}

Applying Lemma \ref{lem2.1}, we next prove the existence of ground states for \eqref{1.6}.
\begin{Proposition} \label{prop2.2}
For every $0\leq\tau\leq1$, the value $\nu_{s,\tau}$ of
\eqref{1.7} is attained at a nonnegative radial decreasing function
$R_\tau$ with $\|R_\tau\|_p=1$, which satisfies
\begin{equation} \label{2.9}
(A_{s,\tau}+1)R_\tau=\nu_{s,\tau}R_\tau^{p-1}\text{ \ in \ }\R^d,
\end{equation}
where the operator $A_{s,\tau}$ is as in \eqref{A}.
Moreover,
\begin{equation} \label{2.10}
 Q_\tau:=\nu_{s,\tau}^{1/(p-2)}R_\tau
\end{equation}
is a positive radial solution of \eqref{1.6}, and the Morse index of $L_{+,\tau}$ defined by \eqref{1.8} is one.
\end{Proposition}

\noindent\textbf{Proof.}
Since $C_c^\infty(\R^d\setminus\{0\})$ is a form core of $\mathfrak h_{s,\tau}$, there exists a minimizing sequence $ \{\phi_n\}\subset C_c^\infty(\R^d\setminus\{0\})$ such that
\begin{equation*}
 \|\phi_n\|_p=1\text{ \ for any \ }n\in\N^+,\ \  \lim_{n\to\infty}\mathfrak h_{s,\tau}[\phi_n]=\nu_{s,\tau}.
\end{equation*}
Let $u_n:=(|\phi_n|)^*$ be the Schwarz rearrangement of $|\phi_n|$. We then obtain from \cite{AlmgrenLieb1989,LiebLoss2001} that for any $u_n\in H^s_{\rm rad}(\R^d)\subset\mathcal Q^s_{\rm rad}$,
\begin{equation*}
\mathfrak h_{s,\tau}[u_n]\leq\mathfrak h_{s,\tau}[\phi_n],\ \ \|u_n\|_p=1.
\end{equation*}
By \eqref{2.7}, there exists a nonnegative radial decreasing function $R_\tau\in D(\mathfrak h_{s,\tau})$ such that up to a subsequence,
\begin{equation*}
 u_n\rightharpoonup R_\tau \text{ \ weakly in \ }D(\mathfrak h_{s,\tau}), \ \
 u_n\to R_\tau\text{ \ strongly in \ }L^p(\R^d)  \text{ \ as \ }n\to\infty.
\end{equation*}
We then derive from the weak lower semicontinuity that
\begin{equation*}
\|R_\tau\|_p=1,\ \ \nu_{s,\tau}\leq\mathfrak h_{s,\tau}[R_\tau]\leq\liminf_{n\to\infty}\mathfrak h_{s,\tau}[u_n]=\nu_{s,\tau}\text{ \ as \  }n\to\infty.
\end{equation*}
The minimizer $R_\tau$ satisfies the Euler--Lagrange equation \eqref{2.9} and is positive by the strong maximum principle (cf. \cite[Theorem~1.1]{JarohsWeth2019}). For $Q_\tau$ in \eqref{2.10}, the constrained second variation implies that
 $\mathfrak{l}_{+,\tau}[\phi]\geq0$ if $\phi\in D(\mathfrak h_{s,\tau})$ satisfies  $\int_{\R^d}Q_\tau^{p-1}\phi dx=0$.
Since this is a codimension-one subspace of the full form domain, the Morse index of $L_{+,\tau}$ is at most one. On the other hand, one can check that
$\mathfrak{l}_{+,\tau}[Q_\tau]=-(p-2)\int_{\R^d}Q_\tau^p dx<0$.
Thus, the Morse index of $L_{+,\tau}$ is exactly one, and the proof of Proposition \ref{prop2.2} is therefore complete.
\qed

 \vspace{0.5em}

 We now employ Proposition \ref{prop2.2} to address the proof of Theorem~\ref{thm:sharp-HGN} (1).
\vspace{0.5em}

\noindent\textbf{Proof of Theorem~\ref{thm:sharp-HGN} (1).}
Let $\theta=d(p-2)/(2sp)$ be as in Theorem~\ref{thm:sharp-HGN}.
For $0\neq u\in\mathcal Q^s$, define $u_\lambda(x):=\lambda^{d/2}u(\lambda x)$ with $\lambda>0$,  so that
\begin{equation}\label{2.13}
 \mathfrak q_{s,1}[u_\lambda]=\lambda^{2s}\mathfrak q_{s,1}[u],\ \ \|u_\lambda\|_2^2=\|u\|_2^2,\ \ \|u_\lambda\|_p^2=\lambda^{2s\theta}\|u\|_p^2.
\end{equation}
Applying the definition of $\nu_{s,1}$ to
$u_\lambda/\|u_\lambda\|_p$, we obtain from \eqref{2.13} that
\begin{equation}\label{2.14}
\nu_{s,1}\|u\|_p^2\leq \lambda^{2s(1-\theta)} \mathfrak q_{s,1}[u]+\lambda^{-2s\theta} \|u\|_2^2.
\end{equation}
Taking the infimum of \eqref{2.14} over $\lambda>0$, we deduce that
$\nu_{s,1}\|u\|_p^2\leq\frac{\mathfrak q_{s,1}[u]^\theta\|u\|_2^{2(1-\theta)}}{\theta^\theta(1-\theta)^{1-\theta}}$,
which thus proves \eqref{1.17}.

Let $R_1$ be the minimizer of $\nu_{s,1}$ obtained in Proposition~\ref{prop2.2}. Since $\|R_1\|_p=1$ and $ \mathfrak h_{s,1}[R_1]=\nu_{s,1}$, consider the $L^p$-preserving dilation $R_{1,\lambda}(x):=\lambda^{d/p}R_1(\lambda x)$. We then have
\begin{equation*}
\mathfrak h_{s,1}[R_{1,\lambda}] =\lambda^{2s(1-\theta)}\mathfrak q_{s,1}[R_1]+\lambda^{-2s\theta}\|R_1\|_2^2.
\end{equation*}
Since $R_1$ is a minimizer under the constraint $\|u\|_p=1$, differentiating  $\mathfrak h_{s,1}[R_{1,\lambda}]$ with respect to $\lambda$  at $\lambda=1$ yields that
$ (1-\theta)\mathfrak q_{s,1}[R_1] = \theta\|R_1\|_2^2$.
Since $ \mathfrak q_{s,1}[R_1]+\|R_1\|_2^2=\nu_{s,1}$, we obtain that
\begin{equation}\label{2.15}
\mathfrak q_{s,1}[R_1]=\theta\nu_{s,1},\ \ \|R_1\|_2^2=(1-\theta)\nu_{s,1}.
\end{equation}
Note from \eqref{2.10} that $\|Q_1\|_p^{p-2}=\nu_{s,1}$. Together with \eqref{2.15}, this implies that $R_1$ attains the identity of   \eqref{1.17} with the best constant \eqref{eq:sharp-HGN-constant-nu}. The proof of Theorem~\ref{thm:sharp-HGN} (1) is therefore complete.
\qed

\begin{Lemma}
\label{lem3.1}
Suppose $d\geq3$, $0<s<1$, $2<p<2_s^*$ and $0<\tau\leq1$. Then every positive solution $Q\in D(\mathfrak h_{s,\tau})$ of \eqref{1.6} belongs to $C^\infty(\R^d\setminus\{0\})$ and satisfies
\begin{equation*}
Q(x)\to+\infty\text{ \ as \ }|x|\to0,\ \ Q(x)\to0\text{ \ as \ }|x|\to\infty.
\end{equation*}
\end{Lemma}

\noindent\textbf{Proof.}
Since $Q\in D(\mathfrak h_{s,\tau})$, we have $Q\in L^2(\R^d)$, and hence
\begin{equation}\label{3.2}
\sup_{x_0\in\R^d}\int_{\R^d}\frac{|Q(y)|}{1+|x_0-y|^{d+2s}}dy\leq C_{d,s}\|Q\|_2.
\end{equation}
Moreover, we obtain from Lemma~\ref{lem2.1} that $Q\in L^q(\R^d)$ holds for $2<q<2_s^*$.
On every ball $B_{2r}(x_0)\Subset\R^d\setminus\{0\}$, equation \eqref{1.6} holds distributionally and $|x|^{-2s}$ is bounded. Since $p<2_s^*$, choose $\max\{2,p-1,d(p-2)/(2s)\}<q<2_s^*$. By applying \cite[Lemma~B.4]{FLS2016} finitely many times, we derive that $Q\in L^\infty_{\rm loc}(\R^d\setminus\{0\})$. Furthermore, let $\Omega\Subset\R^d\setminus\{0\}$ and choose $\eta\in C_c^\infty(\R^d\setminus\{0\})$ satisfying $\eta\equiv1$ near $\overline\Omega$. We have
\begin{equation*}
\begin{aligned}
(-\Delta)^s(\eta Q)=&\bigl(\tau\Cc_{d,s}|x|^{-2s}-1\bigr)Q+Q^{p-1}\\
&+c_{d,s}\int_{\R^d}\frac{(1-\eta(y))Q(y)}{|x-y|^{d+2s}}dy\text{\ \  \ in \  }\mathcal D'(\Omega).
\end{aligned}
\end{equation*}
The integral is smooth on $\Omega$ by \eqref{3.2} and $\operatorname{dist}(\overline\Omega,\operatorname{supp}(1-\eta))>0$. We thus deduce from \cite[Propositions~2.8, 2.9 and 2.22]{Silvestre2007} that $Q\in C^\infty(\R^d\setminus\{0\})$.

Fix sufficiently large $R>0$. Since
$\sup_{|x_0|\geq R}
\big\||x|^{-2s}\big\|_{L^\infty(B_2(x_0))}<\infty$,
by an argument similar to that of \cite[Proposition~B.3]{FLS2016}, we obtain that $\sup_{|x_0|\geq R}\|Q\|_{C^{0,\alpha}(B_1(x_0))}<\infty$ holds for some $\alpha\in(0,1)$.
Together with $Q\in L^2(\R^d)$, this implies that
\begin{equation}\label{3.4}
Q(x)\to0\text{ \ as \ }|x|\to\infty.
\end{equation}
One can check from \cite[Lemma~3.2]{FLSHardy2008} that there exists a unique constant $\delta_\tau\in(0,(d-2s)/2]$ such that
\begin{equation}\label{3.6}
2^{2s}\frac{\Gamma((\delta_\tau+2s)/2)\Gamma((d-\delta_\tau)/2)}{\Gamma(\delta_\tau/2)\Gamma((d-\delta_\tau-2s)/2)}=\tau\Cc_{d,s}.
\end{equation}

Let $\widetilde p_\tau(t,x,y)$ denote the heat kernel of $e^{-tA_{s,\tau}}$. We then derive from \cite[Theorem~1.1]{BogdanEtAl2019} and \eqref{3.6} that for some $c_0>0$,
\begin{equation}\label{radcore:eq-annular-heat-kernel-lower}
\widetilde p_\tau(t,x,y)\geq c_0|x|^{-\delta_\tau}\text{ \ for \  }1\leq t\leq2,\ 0<|x|<\frac12,\ 1<|y|<2.
\end{equation}
Apply the heat-kernel formula to the function $f_n:=\min\{Q^{p-1},n\}\mathbf1_{B_n(0)}$, which converges to $ Q^{p-1}$ a.e. in $\R^d$ as $n\to\infty$. By the monotone convergence theorem, we obtain from Lemma~\ref{lem2.1} that
\begin{equation}\label{eq:nonlinear-heat-representation}
Q(x)=\int_0^\infty e^{-t}\int_{\R^d}\widetilde p_\tau(t,x,y)Q(y)^{p-1} dy dt\text{ \ for \ a.e. \  }x\in\R^d.
\end{equation}
Thus, we derive from \eqref{radcore:eq-annular-heat-kernel-lower} that for $a.e.$  $x\in \R^d$ satisfying $0<|x|<1/2$,
\begin{equation}\label{3.1-pole-lower-bound}
\begin{aligned}
Q(x)\geq&\int_1^2 e^{-t}\int_{1<|y|<2}\widetilde p_\tau(t,x,y)Q(y)^{p-1} dy dt\\
\geq&c_0|x|^{-\delta_\tau}\int_1^2 e^{-t} dt\int_{1<|y|<2}Q(y)^{p-1} dy\geq C|x|^{-\delta_\tau}.
\end{aligned}
\end{equation}
Since $Q$ is continuous on $\R^d\setminus\{0\}$, the lower bound $Q(x)\geq C|x|^{-\delta_\tau}$ holds for $x\in \R^d$ satisfying $0<|x|<1/2$.
This completes the proof of Lemma~\ref{lem3.1} in view of \eqref{3.4} and \eqref{3.1-pole-lower-bound}.
\qed

\begin{Lemma}
\label{lem2.4}
Suppose $d\geq3$, $0<s<1$, and $2<p<2_s^*$. For every $0<\tau\leq1$, every positive solution $Q\in D(\mathfrak h_{s,\tau})$ of \eqref{1.6} is radially symmetric about the origin and strictly decreasing in $r=|x|$. At $\tau=0$, every positive solution of \eqref{1.6} in $H^s(\R^d)$ is, up to translation, radial and strictly decreasing.
\end{Lemma}

\noindent\textbf{Proof.}
The case $\tau=0$ follows directly from \cite[Proposition~3.1]{FLS2016}. We then use a moving-plane argument (cf.~\cite{ChenLiLi2017,MaZhao2010}) to consider the case $0<\tau\leq1$. Set for $\lambda>0$,
\begin{equation*}
\Sigma_\lambda:=\{x\in\R^d:x_1>\lambda\},\ \ w_\lambda(x):=Q(x^\lambda)-Q(x),\ \ a_\lambda:=2\lambda e_1,
\end{equation*}
where $x^\lambda:=(2\lambda-x_1,x')$. It then follows from \eqref{3.1-pole-lower-bound} that
\begin{equation}\label{eq:2.21}
w_\lambda(x)\to+\infty\text{ \ as \ }x\to a_\lambda.
\end{equation}
Moreover, we obtain from \eqref{1.6} that
\begin{equation}\label{eq:2.22}
(-\Delta)^sw_\lambda+c_\lambda(x)w_\lambda=\tau\Cc_{d,s}\bigl(|x^\lambda|^{-2s}-|x|^{-2s}\bigr)Q(x^\lambda)\geq0\text{ \ in \ }\Sigma_\lambda\setminus\{a_\lambda\},
\end{equation}
where
\begin{equation}\label{eq:2.23}
c_\lambda(x):=1-\tau\Cc_{d,s}|x|^{-2s}-b_\lambda(x),
\end{equation}
and
\begin{equation*}
b_\lambda(x):=\begin{cases}
\displaystyle\frac{Q(x^\lambda)^{p-1}-Q(x)^{p-1}}{Q(x^\lambda)-Q(x)},& \ \ \text{if}~~Q(x^\lambda)\neq Q(x); \\
(p-1)Q(x)^{p-2},& \ \ \text{if}~~ Q(x^\lambda)=Q(x).
\end{cases}
\end{equation*}
Denoting
\begin{equation*}
D_\lambda:=\big\{x\in\Sigma_\lambda:~w_\lambda(x)<0\big\},
\end{equation*}
we have
\begin{equation}\label{eq:negative}
0\leq b_\lambda(x)\leq(p-1)Q(x)^{p-2}\ \ \text{ \ in }D_\lambda.
\end{equation}
Suppose $D_\lambda\neq\varnothing$, then we have
$m:=\inf_{x\in\Sigma_\lambda}w_\lambda(x)<0$.
Moreover, $w_\lambda=0$ on $\partial\Sigma_\lambda$ and $w_\lambda(x)\to0$ as $|x|\to\infty$. Therefore, we obtain from \eqref{eq:2.21} that $m$ is attained at some point $z\in\Sigma_\lambda$. In particular, we have $w_\lambda(z)=m<0$, and hence $z\in D_\lambda$.
By the antisymmetry of $w_\lambda$, one can check that
\begin{equation}\label{l0}
\begin{aligned}
(-\Delta)^sw_\lambda(z)=&c_{d,s}\operatorname{p.v.}\int_{\Sigma_\lambda}(m-w_\lambda(y))\left(\frac1{|z-y|^{d+2s}}-\frac1{|z-y^\lambda|^{d+2s}}\right)dy\\
&+2c_{d,s}m\int_{\Sigma_\lambda}\frac{dy}{|z-y^\lambda|^{d+2s}}
\leq \gamma_{d,s}(z_1-\lambda)^{-2s}m<0,
\end{aligned}
\end{equation}
where $\gamma_{d,s}>0$ depends only on $d$ and $s$.

Since $|x|\geq\lambda$ in $\Sigma_\lambda$ and $Q(x)\to0$ as $|x|\to\infty$, we obtain from \eqref{eq:2.23} and \eqref{eq:negative} that  for sufficiently large $\lambda>0$, $c_\lambda\geq1/2$ holds in $D_\lambda$. We then conclude from \eqref{eq:2.22} and \eqref{l0} that  for all sufficiently large $\lambda>0$,
\begin{equation*}
w_\lambda\geq0\text{ \ in \ }\Sigma_\lambda.
\end{equation*}
Denote
\begin{equation*}
\lambda_0:=\inf\left\{\lambda>0:~~w_\mu\geq0\text{ in }\Sigma_\mu\text{ for every }\mu\geq\lambda\right\},
\end{equation*}
and suppose $\lambda_0>0$. One can check that $w_{\lambda_0}\geq0$ in $\Sigma_{\lambda_0}\setminus\{a_{\lambda_0}\}$ by the continuity of $Q$ away from the origin. If $w_{\lambda_0}(z)=0$ holds for some $z\in\Sigma_{\lambda_0}\setminus\{a_{\lambda_0}\}$, then we obtain from \eqref{eq:2.22} and the antisymmetry that
\begin{equation*}
\begin{aligned}
0\geq&c_{d,s}\int_{\Sigma_{\lambda_0}}w_{\lambda_0}(y)\left(\frac1{|z-y^{\lambda_0}|^{d+2s}}-\frac1{|z-y|^{d+2s}}\right)dy\\
=&(-\Delta)^sw_{\lambda_0}(z)=\tau\Cc_{d,s}\bigl(|z^{\lambda_0}|^{-2s}-|z|^{-2s}\bigr)Q(z^{\lambda_0})>0,
\end{aligned}
\end{equation*}
a contradiction. We therefore deduce that
\begin{equation}\label{eq:reflec}
w_{\lambda_0}>0\text{ \ in \ }\Sigma_{\lambda_0}\setminus\{a_{\lambda_0}\}.
\end{equation}

Applying the decay of $Q$ as $|x|\rightarrow \infty$, we obtain from \eqref{eq:2.23} and \eqref{eq:negative} that there exists $R>3\lambda_0$ such that for $\lambda>0$,
\begin{equation}\label{eq:farpositive}
c_\lambda(x)\geq0\text{ \ in \ } D_\lambda\cap\{|x|>R\}.
\end{equation}
We also derive from \eqref{eq:2.23} and \eqref{eq:negative} that for $3\lambda_0/4\leq\lambda\leq\lambda_0$ and $x\in D_\lambda\cap\{|x|\leq R\}$,
\begin{equation}\label{02}
c_\lambda(x)\geq1-\tau\Cc_{d,s}\left(\frac{\lambda_0}{2}\right)^{-2s}-(p-1)\sup_{\lambda_0/2\leq|y|\leq R}Q(y)^{p-2}\geq-M,
\end{equation}
where $M>0$ is independent of $\lambda$.
Choose $0<\rho<\min\{1/4,\lambda_0/4\}$ such that
$$C\rho^{-\delta_\tau}>\sup_{\lambda_0\leq|y|\leq3\lambda_0}Q(y),$$
where $C>0$ and $\delta_\tau>0$ are as in \eqref{3.1-pole-lower-bound}. We also have $|x^\lambda|=|x-a_\lambda|<\rho$ for $3\lambda_0/4\leq\lambda\leq\lambda_0$ and $0<|x-a_\lambda|<\rho$. Moreover, since $|a_\lambda|=2\lambda$ and $\rho<\lambda_0/4$, the triangle inequality yields that $\lambda_0<|x|<3\lambda_0$. Therefore, we get from \eqref{3.1-pole-lower-bound} that
\begin{equation}\label{eq:unip}
w_\lambda(x)\geq C\rho^{-\delta_\tau}-\sup_{\lambda_0\leq|y|\leq3\lambda_0}Q(y)>0.
\end{equation}

Choose $0<\varepsilon<\lambda_0/4$ such that $\gamma_{d,s}(2\varepsilon)^{-2s}>M$, and set
\begin{equation*}
K_\varepsilon:=\left\{x\in\R^d:~|x|\leq R,\ \ x_1\geq\lambda_0+\varepsilon\right\}\setminus B_{\f{\rho}{2}}(a_{\lambda_0}).
\end{equation*}
We then obtain from \eqref{eq:reflec} that $\min_{K_\varepsilon}w_{\lambda_0}>0$. If $|a_\lambda-a_{\lambda_0}|<\rho/4$, then $|x^\lambda|\geq\rho/4$ holds on $K_\varepsilon$. Thus, $w_\lambda\to w_{\lambda_0}$ uniformly on $K_\varepsilon$ as $\lambda\to\lambda_0$. Taking $\lambda\in[3\lambda_0/4,\lambda_0)$ sufficiently close to $\lambda_0$, we then have
\begin{equation*}
\lambda_0-\lambda<\varepsilon,\ \ |a_\lambda-a_{\lambda_0}|<\frac{\rho}{4},\ \ \min_{K_\varepsilon}w_\lambda>0.
\end{equation*}
If $D_\lambda\neq\varnothing$, let $z$ be a negative minimum point as in \eqref{l0}. We then deduce from \eqref{eq:2.22}, \eqref{l0}, \eqref{eq:farpositive} and \eqref{eq:unip} that $|z|\leq R$ and $z\notin B_\rho(a_\lambda)$. Since $B_{\rho/2}(a_{\lambda_0})\subset B_\rho(a_\lambda)$ and $w_\lambda>0$ on $K_\varepsilon$, one can check that
\begin{equation*}
z_1<\lambda_0+\varepsilon,\ \ 0<z_1-\lambda<2\varepsilon,
\end{equation*}
where $z=(z_1,\cdots,z_d)\in \R^d$.
We then deduce from \eqref{eq:2.22}, \eqref{l0} and \eqref{02}  that
\begin{equation*}
0\leq(-\Delta)^sw_\lambda(z)+c_\lambda(z)m\leq m[\gamma_{d,s}(2\varepsilon)^{-2s}-M]<0,
\end{equation*}
a contradiction. Thus, $D_\lambda=\varnothing$ holds for all $\lambda<\lambda_0$ sufficiently close to $\lambda_0$, which contradicts the definition of $\lambda_0$.

We conclude from above that $\lambda_0=0$. The argument of \eqref{eq:reflec} then yields that  for every $\lambda>0$, $w_\lambda>0$ holds in $\Sigma_\lambda\setminus\{a_\lambda\}$. Repeating the argument in the opposite direction and then in every direction $e\in\mathbb S^{d-1}$, we derive that $Q$ is symmetric about every hyperplane through the origin. For $0<r_1<r_2$, the strict reflection across $x\cdot e=(r_1+r_2)/2$ implies that $Q(r_1e)>Q(r_2e)$. This therefore completes the proof of Lemma~\ref{lem2.4}.
\qed

 \vspace{0.5em}
The next lemma is concerned with the exact asymptotics of positive solutions for \eqref{1.6} as $|x|\to0$.
\begin{Lemma}
\label{prop:critical-pole-asymptotics}
Suppose $d\geq3$, $0<s<1$, and $2<p<2_s^*$.  If
$0<Q\in\mathcal Q^s$ solves \eqref{1.6} with $\tau=1$, then there exists
$\ell_Q\in(0,\infty)$ such that
\begin{equation}\label{eq:exaasy}
\lim_{|x|\to0}|x|^{(d-2s)/2}Q(x)=\ell_Q.
\end{equation}
\end{Lemma}

\noindent\textbf{Proof.}
Set
\begin{equation*}
\alpha:=2s,\quad \delta:=\frac{d-\alpha}{2}=\frac{d-2s}{2},\quad h(x):=|x|^{-\delta}.
\end{equation*}
By Lemma~\ref{lem2.4}, $Q$ is radial and strictly decreasing in $|x|>0$. Denote $m:=d-d/p$. We then have
\begin{equation}\label{eq:nonbond}
Q(r)^{p-1}\leq\Big(\frac1{|B_r(0)|}\int_{B_r(0)}Q^p dx\Big)^{(p-1)/p}\leq C_Qr^{-m}\text{ \ for \  }r>0.
\end{equation}
Since $2<p<2_s^*=2d/(d-\alpha)$, we have
\begin{equation}\label{eq:exponent}
\delta+m<d, \ \ -1<\frac{\delta-m}{\alpha}<0.
\end{equation}

Let $\widetilde p_1(t,x,y)$ denote the heat kernel of the semigroup
$e^{-tA_{s,1}}$.  It then follows from \cite[Theorem~1.1]{BogdanEtAl2019} that
\begin{equation} \label{eq:kerupp}
\widetilde p_1(t,x,y)\leq C t^{-d/\alpha}\left(1\wedge\big|t^{-1/\alpha}(x-y)\big|^{-d-\alpha}\right)\left(1+t^{\delta/\alpha}|x|^{-\delta}\right)\left(1+t^{\delta/\alpha}|y|^{-\delta}\right),
\end{equation}
where and below we denote $a\wedge b:=\min\{a,b\}$.
Moreover, the symmetry of $\widetilde p_1$ and
\cite[Theorem~1.1]{BogdanJakubowskiKimPilarczyk2023} imply that for every
$t>0$ and $y\neq0$,
\begin{equation} \label{eq:kerlim}
 \Psi_t(y):=\lim_{|x|\to0}\frac{\widetilde p_1(t,x,y)}{h(x)}\in(0,\infty).
\end{equation}
Define $F_t(x):=\int_{\R^d}\widetilde p_1(t,x,y)Q(y)^{p-1}dy/h(x)$ and set $\xi:=t^{-1/\alpha}x$. Using $y=t^{1/\alpha}z$, we then obtain from \eqref{eq:nonbond} and \eqref{eq:kerupp} that
\begin{equation*}
\begin{aligned}
F_t(x)\leq&C t^{(\delta-m)/\alpha}(1+|\xi|^\delta)\int_{\R^d}\left(1\wedge|\xi-z|^{-d-\alpha}\right)(1+|z|^{-\delta})|z|^{-m}dz\\
\leq&C t^{(\delta-m)/\alpha}(1+|\xi|)^{\delta-m},
\end{aligned}
\end{equation*}
where the last inequality follows by splitting the integral into $|z|\leq2$ and $|z|>2$ for $|\xi|\leq1$, and into $|z-\xi|\leq|\xi|/2$, $|z|\leq|\xi|/2$ and their complement for $|\xi|>1$. Since $\delta-m<0$, we have
\begin{equation}\label{eq:larget}
e^{-t}F_t(x)\leq Ce^{-t}t^{(\delta-m)/\alpha}\text{ \ for \  }t>0,\ x\neq0.
\end{equation}
Applying dominated convergence theorem, it yields from \eqref{eq:kerupp} and \eqref{eq:kerlim} that for every fixed $t>0$,
\begin{equation*}
F_t(x)\to\int_{\R^d}\Psi_t(y)Q(y)^{p-1}dy\text{ \ as \  }|x|\to0.
\end{equation*}
The right-hand side of \eqref{eq:larget} is integrable on $(0,\infty)$ in view of \eqref{eq:exponent}. Using dominated convergence theorem and the continuity of $Q/h$ on $\R^d\setminus\{0\}$, we then obtain from \eqref{eq:nonlinear-heat-representation} and \eqref{eq:larget} that
\begin{equation*}
\lim_{|x|\to0}\frac{Q(x)}{h(x)}=\int_0^\infty e^{-t}\int_{\R^d}\Psi_t(y)Q(y)^{p-1}dy dt=:\ell_Q\text{ \ as \  }|x|\to0.
\end{equation*}
Finally, we deduce from \eqref{eq:exponent} and \eqref{eq:larget} that $\ell_Q<\infty$, while it implies from \eqref{eq:kerlim} and $Q>0$ that $\ell_Q>0$. This therefore completes the proof of Lemma~\ref{prop:critical-pole-asymptotics}.
\qed

\vspace{0.5em}

\noindent\textbf{Proof of Theorem~\ref{thm1.1}.}
Proposition~\ref{prop2.2} gives the existence of ground states for \eqref{1.6}. Lemmas~\ref{lem3.1} and \ref{lem2.4} give the smoothness away from $x=0$, the decaying property at infinity, radial symmetry, and strict radial monotonicity of positive weak solutions for \eqref{1.6}. Lemma \ref{prop:critical-pole-asymptotics} establishes the exact asymptotic formula \eqref{eq:exaasy}, which therefore proves Theorem~\ref{thm1.1}.
\qed

\section*{Declaration of AI-assisted tools}
During the preparation of this paper, OpenAI models were used to assist with exploratory mathematical discussions and analyze the orthogonality condition \eqref{eq:2moment}. The authors take full responsibility for mathematical contents and the final version of the manuscript.

\end{document}